\documentclass[final]{siamltex}

\usepackage{amssymb}
\usepackage{amsmath}
\usepackage{graphicx}
\usepackage{subfigure}
\usepackage{color}
\usepackage{psfrag}
\usepackage{url}

\newcommand{\mb}[1]{\mbox{\boldmath$#1$}}

\begin{document}

\title{Angular Synchronization by Eigenvectors and Semidefinite Programming}
\author{
A.~Singer%
\thanks{Department of Mathematics and PACM, Princeton University, Fine Hall, Washington Road, Princeton NJ 08544-1000 USA, Email: amits@math.princeton.edu}}

\date{}
\maketitle

\begin{abstract}
The angular synchronization problem is to obtain an accurate estimation (up to a constant additive phase) for a set of unknown angles $\theta_1,\ldots,\theta_n$ from $m$ noisy measurements of their offsets $\theta_i-\theta_j \mod 2\pi$. Of particular interest is angle recovery in the presence of many outlier measurements that are uniformly distributed in $[0,2\pi)$ and carry no information on the true offsets. We introduce an efficient recovery algorithm for the unknown angles from the top eigenvector of a specially designed Hermitian matrix. The eigenvector method is extremely stable and succeeds even when the number of outliers is exceedingly large. For example, we successfully estimate $n=400$ angles from a full set of $m={400 \choose 2}$ offset measurements of which $90\%$ are outliers in less than a second on a commercial laptop. The performance of the method is analyzed using random matrix theory and information theory. We discuss the relation of the synchronization problem to the combinatorial optimization problem \textsc{Max-2-Lin mod} $L$ and present a semidefinite relaxation for angle recovery, drawing similarities with the Goemans-Williamson algorithm for finding the maximum cut in a weighted graph. We present extensions of the eigenvector method to other synchronization problems that involve different group structures and their applications, such as the time synchronization problem in distributed networks and the surface reconstruction problems in computer vision and optics.
\end{abstract}

\section{Introduction}
\label{sec:intro}

The angular synchronization problem is to estimate $n$ unknown angles $\theta_1,\ldots,\theta_n \in [0,2\pi)$ from $m$ noisy measurements $\delta_{ij}$ of their offsets $\theta_i-\theta_j \mod 2\pi$. In general, only a subset of all possible $n \choose 2$ offsets are measured. The set $E$ of pairs $\{i,j\}$ for which offset measurements exist can be realized as the edge set of a graph $G=(V,E)$ with vertices corresponding to angles and edges corresponding to measurements.

When all offset measurements are exact with zero measurement error, it is possible to solve the angular synchronization problem iff the graph $G$ is connected. Indeed, if $G$ is connected then it contains a spanning tree and all angles are sequentially determined by traversing the tree while summing the offsets modulo $2\pi$. The angles are uniquely determined up to an additive phase, e.g., the angle of the root. On the other hand, if $G$ is disconnected then it is impossible to determine the offset between angles that belong to disjoint components of the graph.

Sequential algorithms that integrate the measured offsets over a particular spanning tree of the graph are very sensitive to measurement errors, due to accumulation of the errors. It is therefore desirable to integrate all offset measurements in a globally consistent way. The need for such a globally consistent integration method comes up in a variety of applications. One such application is the time synchronization of distributed networks \cite{GiridharKumar,Karp}, where clocks measure noisy time offsets $t_i-t_j$ from which the determination of $t_1,\ldots,t_n\in \mathbb{R}$ is required. Other applications include the surface reconstruction problems in computer vision \cite{FrankotChellappa,AgrawalWhatRange} and optics \cite{Koby2001}, where the surface is to be reconstructed from noisy measurements of the gradient to the surface and the graph of measurements is typically the two-dimensional regular grid. The most common approach in the above mentioned applications for a self consistent global integration is the least squares approach. The least squares solution is most suitable when the offset measurements have a small Gaussian additive error. The least squares solution can be efficiently computed and also mathematically analyzed in terms of the Laplacian of the underlying measurement graph.

There are many possible models for the measurement errors, and we are mainly interested in models that allow many outliers. An outlier is an offset measurement that has a uniform distribution on $[0,2\pi)$ regardless of the true value for the offset. In addition to outliers that carry no information on the true angle values, there also exist of course good measurements whose errors are relatively small. We have no a-priori knowledge, however, which measurements are good and which are bad (outliers).

In our model, the edges of $E$ can be split into a set of good edges $E_{good}$ and a set of bad edges $E_{bad}$, of sizes $m_{good}$ and $m_{bad}$ respectively (with $m = |E| =  m_{good}+m_{bad}$), such that
\begin{equation}
\label{model}
\begin{array}{cc}
\delta_{ij} = \theta_i - \theta_j & \mbox{for } \{i,j\} \in E_{good} \\
\delta_{ij} \sim Uniform\left([0,2\pi)\right) & \mbox{for } \{i,j\} \in E_{bad}
\end{array}.
\end{equation}
Perhaps it would be more realistic to allow a small discretization error for the good offsets, for example, by letting them have the wrapped normal distribution on the circle with mean $\theta_i-\theta_j$ and variance $\sigma^2$ (where $\sigma$ is a typical discretization error).  This discretization error can be incorporated into the mathematical analysis of Section \ref{sec:random} with a little extra difficulty. However, the effect of the discretization error is negligible compared to that of the outliers, so we choose to ignore it in order to make the presentation as simple as possible.

It is trivial to find a solution to (\ref{model}) if some oracle whispers to our ears which equations are good and which are bad (in fact, all we need in that case is that $E_{good}$ contains a spanning tree of $G$). In reality, we have to be able to tell the good from the bad on our own.

The overdetermined system of linear equations (modulo $2\pi$)
\begin{equation}
\label{theta_E}
\theta_i - \theta_j = \delta_{ij} \mod 2\pi, \quad \mbox{for } \{i,j\}\in E
\end{equation}
can be solved by the method of least squares as follows. Introducing the complex-valued variables $z_i=e^{\imath \theta_i}$, the system (\ref{theta_E}) is equivalent to
\begin{equation}
\label{eq:lsqr}
z_i - e^{\imath \delta_{ij}}z_j = 0, \quad \mbox{for } \{i,j\}\in E,
\end{equation}
which is an overdetermined system of homogeneous linear equations over $\mathbb{C}$. To prevent the solution from collapsing to the trivial solution $z_1=z_2=\cdots = z_n = 0$, we set $z_1=1$ (recall that the angles are determined up to a global additive phase, so we may choose $\theta_1=0$), and look for the solution $z_2,\ldots,z_n$ of (\ref{eq:lsqr}) with minimal $\ell_2$-norm residual.
However, it is expected that the sum of squares errors would be overwhelmingly dominated by outlier equations, making least squares least favorable to succeed if the proportion of bad equations is large (see numerical results involving least squares in Table \ref{tab:sdp-eig2}). We therefore seek for a solution method which is more robust to outliers.

Maximum likelihood is an obvious step in that direction. The maximum likelihood solution to (\ref{model}) is simply the set of angles $\theta_1,\ldots,\theta_n$ that satisfies as many equations of (\ref{theta_E}) as possible. We may therefore define the self consistency error (SCE) of $\theta_1,\ldots,\theta_n$ as the number of equations not being satisfied
\begin{equation}
\label{SCE}
SCE(\theta_1,\ldots,\theta_n) = \#\{\{i,j \}\in E : \theta_i - \theta_j \neq \delta_{ij} \mod 2\pi\}.
\end{equation}
As even the good equations contain some error (due to angular discretization and noise), a more suitable self consistency error is $SCE_f$ that incorporates some penalty function $f$
\begin{equation}
\label{SCEf}
SCE_f(\theta_1,\ldots,\theta_n) = \sum_{\{ i,j \}\in E} f(\theta_i-\theta_j - \delta_{ij}),
\end{equation}
where $f:[0,2\pi)\to \mathbb{R}$ is a smooth periodic function with $f(0)=0$ and $f(\theta)=1$ for $|\theta|>\theta_0$, where $\theta_0$ is the allowed discretization error. The minimization of (\ref{SCEf}) is equivalent to maximizing the log likelihood with a different probabilistic error model.

The maximum likelihood approach suffers from a major drawback though. It is virtually impossible to find the global minimizer $\theta_1,\ldots,\theta_n$ when dealing with large scale problems ($n\gg 1$), because the minimization of either (\ref{SCE}) or (\ref{SCEf}) is a non-convex optimization problem in a huge parameter space.
It is like finding a needle in a haystack. 

In this paper we take a different approach and introduce two different estimators for the angles. The first estimator is based on an eigenvector computation while the second estimator is based on a semidefinite program (SDP) \cite{Boyd1996}. Our eigenvector estimator $\hat{\theta}_1,\ldots,\hat{\theta}_n$ is obtained by the following two-step recipe . In the first step, we construct an $n\times n$ complex-valued matrix $H$ whose entries are
\begin{equation}
\label{introH}
H_{ij} = \left\{\begin{array}{cc} e^{\imath \delta_{ij}} & \{i,j\}\in E \\
0 & \{i,j\} \not \in E
\end{array} \right.,
\end{equation}
where $\imath = \sqrt{-1}$. The matrix $H$ is Hermitian, i.e. $H_{ij} = \bar{H}_{ji}$, because the offsets are skew-symmetric $\delta_{ij}=-\delta_{ji} \mod 2\pi$. As $H$ is Hermitian, its eigenvalues are real. The second step is to compute the top eigenvector $v_1$ of $H$ with maximal eigenvalue, and to define the estimator in terms of this top eigenvector as
\begin{equation}
\label{eig}
e^{\imath \hat{\theta}_i} = \frac{v_1(i)}{|v_1(i)|},\quad i=1,\ldots,n.
\end{equation}
The philosophy leading to the eigenvector method is explained in  Section \ref{sec:eigenvector}.

The second estimator is based on the following SDP
\begin{align}
& \max_{\Theta \in \mathbb{C}^{n\times n}} \operatorname{trace}(\bar{H}\Theta) \label{sdp1a}\\
s.t. & \Theta \succeq 0 \label{sdp1b}\\
& \Theta_{ii}=1 \quad i=1,2,\ldots,n,\label{sdp1c}
\end{align}
where $\Theta \succeq 0$ is a shorthand notation for $\Theta$ being a Hermitian semidefinite positive matrix. The only difference between this SDP and the Goemans-Williamson algorithm for finding the maximum cut in a weighted graph \cite{GoemansWilliamson} is that the maximization is taken over all semidefinite positive Hermitian matrices with complex-valued entries rather than just the real-valued symmetric matrices. The SDP-based estimator $\hat{\theta}_1,\ldots,\hat{\theta}_n$ is derived from the normalized top eigenvector $v_1$ of $\Theta$ by the same rounding procedure (\ref{eig}).
Our numerical experiments show that the accuracy of the eigenvector method and the SDP method are comparable. Since the eigenvector method is much faster, we prefer using it for large scale problems. The eigenvector method is also numerically appealing, because in the useful case the spectral gap is large, rendering the simple power method an efficient and numerically stable way of computing the top eigenvector. The SDP method is summarized in Section \ref{sec:sdp}.

In Section \ref{sec:random} we use random matrix theory to analyze the eigenvector method for two different measurement graphs: the complete graph and ``small-world" graphs \cite{small_world}. Our analysis shows that the top eigenvector of $H$ in the complete graph case has a non-trivial correlation with the vector of true angles as soon as the proportion $p$ of good offset measurements becomes greater than $\frac{1}{\sqrt{n}}$. In particular, the correlation goes to 1 as ${np^2}\to \infty$, meaning a successful recovery of the angles. Our numerical simulations confirm these results and demonstrate the robustness of the estimator (\ref{eig}) to outliers.

In Section \ref{sec:information-theory} we prove that the eigenvector method is asymptotically nearly optimal in the sense that it achieves the information theoretic Shannon bound up to a multiplicative factor that depends only on the discretization error of the measurements $2\pi/L$, but not on $m$ and $n$. In other words, no method whatsoever can accurately estimate the angles if the proportion of good measurements is $o(\sqrt{\frac{n}{m}})$.
The connection between the angular synchronization problem and \textsc{Max-2-Lin mod} $L$ \cite{Hastad} is explored in Section \ref{sec:max2lin}. Finally, Section \ref{sec:sync} is a summary and discussion of further applications of the eigenvector method to other synchronization problems over different groups.

\section{The Eigenvector Method}
\label{sec:eigenvector}

Our approach to finding the self consistent solution for $\theta_1,\ldots,\theta_n$ starts with forming the following $n \times n$ matrix $H$
\begin{equation}
\label{H}
H_{ij} = \left\{\begin{array}{cc} e^{\imath \delta_{ij}} & \{i,j\}\in E \\
0 & \{i,j \} \not\in E \end{array} \right.,
\end{equation}
where $\imath = \sqrt{-1}$. Since
\begin{equation}
\label{hermitian}
\delta_{ij} = -\delta_{ij},\quad \mbox{for all } i,j=1,\ldots,n,
\end{equation}
it follows that $H_{ij} = \bar{H}_{ji}$, where for any complex number $z=a+\imath b$ we denote by $\bar{z}=a-\imath b$ its complex conjugate. In other words, the matrix $H$ is Hermitian, i.e., $H^*=H$.

Next, we consider the maximization problem
\begin{equation}
\label{max}
\max_{\theta_1,\ldots,\theta_n \in [0,2\pi)}\sum_{i,j=1}^n e^{-\imath\theta_i} H_{ij} e^{\imath\theta_j},
\end{equation}
and explain the philosophy behind it. For the correct set of angles $\theta_1,\ldots,\theta_n$, each good edge contributes $$e^{-\imath \theta_i} e^{\imath (\theta_i -\theta_j)} e^{\imath \theta_j}=1$$ to the sum in (\ref{max}).  The total contribution of the good edges is just the sum of ones, piling up to be exactly the total number of good edges $m_{good}$. On the other hand, the contribution of each bad edge will be uniformly distributed on the unit circle in the complex plane. Adding up the terms due to bad edges can be thought of as a discrete planar random walk where each bad edge corresponds to a unit size step at a uniformly random direction. These random steps mostly cancel out each other, such that the total contribution of the $m_{bad}$ edges is only $O(\sqrt{m_{bad}})$. It follows that the objective function in (\ref{max}) has the desired property of diminishing the contribution of the bad edges by a square root relative to the linear contribution of the good edges.

Still, the maximization problem (\ref{max}) is a non-convex maximization problem which is quite difficult to solve in practice. We therefore introduce the following relaxation of the problem
\begin{equation}
\label{max2}
\max_{\begin{array}{c}z_1,\ldots,z_n \in \mathbb{C} \\ \sum_{i=1}^n |z_i|^2 = n \end{array}}\sum_{i,j=1}^n z_i^* H_{ij} z_j.
\end{equation}
That is, we replace the previous $n$ individual constraints for each of the variables $z_i=e^{\imath \theta_i}$ to have a unit magnitude, by a single and much weaker constraint, requiring the sum of squared magnitudes to be $n$. The maximization problem (\ref{max2}) is that of a quadratic form whose solution is simply given by the top eigenvector of the Hermitian matrix $H$. Indeed, the spectral theorem implies that the eigenvectors $v_1,v_2,\ldots,v_n$ of $H$ form an orthonormal basis for $\mathbb{C}^n$ with corresponding real eigenvalues $\lambda_1 \geq \lambda_2 \geq \ldots \geq \lambda_n$ satisfying $Hv_i = \lambda_i v_i$. Rewriting the constrained maximization problem (\ref{max2}) as
\begin{equation}
\max_{\|z\|^2 = n} z^* H z,
\end{equation}
it becomes clear that the maximizer $z$ is given by $z=v_1$, where $v_1$ is the normalized top eigenvector satisfying $Hv_1 = \lambda_1 v_1$ and $\|v_1\|^2=n$, with $\lambda_1$ being the largest eigenvalue. The components of the eigenvector $v_1$ are not necessarily of unit magnitude, so we normalize them and define the estimated angles by
\begin{equation}
\label{eig2}
e^{\imath \hat{\theta}_i} = \frac{v_1(i)}{|v_1(i)|},\quad \mbox{for } i=1,\ldots,n
\end{equation}
(see also equation (\ref{eig})).

The top eigenvector can be efficiently computed by the power iteration method that starts from a randomly chosen vector $b_0$ and iterates $b_{n+1} = \frac{Hb_n}{\|Hb_n \|}$. Each iteration requires just a matrix-vector multiplication that takes $O(n^2)$ operations for dense matrices, but only $O(m)$ operations for sparse matrices, where $m=|E|$ is the number of non-zero entries of $H$ corresponding to edges in the graph. The number of iterations required by the power method decreases with the spectral gap that indeed exists and is analyzed in detail in Section \ref{sec:random}.

Note that cycles in the graph of good edges lead to consistency relations between the offset measurements. For example, if the three edges $\{i,j\},\{j,k\},\{k,i\}$ are a triangle of good edges, then the corresponding offset angles $\delta_{ij}$, $\delta_{jk}$ and $\delta_{ki}$ must satisfy
\begin{equation}
\label{consistent}
\delta_{ij} + \delta_{jk} + \delta_{ki} = 0 \mod 2\pi,
\end{equation}
because
$$\delta_{ij} + \delta_{jk} + \delta_{ki} = \theta_i - \theta_j + \theta_j - \theta_k + \theta_k - \theta_i = 0 \mod 2\pi.$$
A closer look into the power iteration method reveals that multiplying the matrix $H$ by itself integrates the information in the consistency relation of triplets, while higher order iterations exploit consistency relations of longer cycles.
Indeed,
\begin{eqnarray}
H^2_{ij} &=& \sum_{k=1}^n H_{ik} H_{kj} = \sum_{k: \{i,k \},\{j,k\} \in E} e^{\imath \delta_{ik}} e^{\imath \delta_{kj}} = \sum_{k: \{i,k \},\{j,k\} \in E} e^{-\imath (\delta_{jk}+\delta_{ki})} \label{H1}\\
&=& \#\left\{k: \{i,k \} \mbox{ and } \{j,k\} \in E_{good}\right\} e^{\imath (\theta_i-\theta_j)} \label{H2}\\ && + \sum_{k: \{i,k \} \mbox{ or } \{j,k\} \in E_{bad}} e^{-\imath (\delta_{jk}+\delta_{ki})},\nonumber
\end{eqnarray}
where we employed (\ref{hermitian}) in (\ref{H1}), and (\ref{consistent}) in (\ref{H2}).
The top eigenvector therefore integrates the consistency relations of all cycles.

\section{The semidefinite program approach}
\label{sec:sdp}
A different natural relaxation of the optimization problem (\ref{max}) is using SDP. Indeed, the objective function in (\ref{max}) can be written as
\begin{equation}
\sum_{i,j=1}^n e^{-\imath\theta_i} H_{ij} e^{\imath\theta_j} = \operatorname{trace}(\bar{H}\Theta),
\end{equation}
where $\Theta$ is the $n\times n$ complex-valued rank-one Hermitian matrix
\begin{equation}
\label{rank-one}
\Theta_{ij} = e^{\imath (\theta_i - \theta_j)}.
\end{equation}
Note that $\Theta$ has ones on its diagonal
\begin{equation}
\Theta_{ii} = 1,\quad i=1,2,\ldots,n.
\end{equation}
Except for the non-convex rank-one constraint implied by (\ref{rank-one}), all other constraints are convex and lead to the natural SDP relaxation (\ref{sdp1a})-(\ref{sdp1c}). This program is almost identical to the Goemans-Williamson SDP for finding the maximum cut in a weighted graph. The only difference is that here we maximize over all possible complex-valued Hermitian matrices, not just the symmetric real matrices. The SDP-based estimator corresponding to (\ref{sdp1a})-(\ref{sdp1c}) is then obtained from the best rank-one approximation of the optimal matrix $\Theta$ using the Cholesky decomposition.

The SDP method may seem favorable to the eigenvector method as it explicitly imposes the unit magnitude constraint for $e^{\imath \theta_i}$. Our numerical experiments show that the two methods give similar results (see Table \ref{tab:sdp-eig2}). Since the eigenvector method is much faster, it is also the method of choice for large scale problems.

\section{Connections with random matrix theory and spectral graph theory}
\label{sec:random}

In this section we analyze the eigenvector method using tools from random matrix theory and spectral graph theory.

\subsection{Analysis of the complete graph angular synchronization problem}
We first consider the angular synchronization problem in which all ${n \choose 2}$ angle offsets are given, so that the corresponding graph is the complete graph $K_n$ of $n$ vertices. We also assume that the probability for each edge to be good is $p$, independently of all other edges.
This probabilistic model for the graph of good edges is known as the Erd\H{o}s-R\'enyi random graph $G(n,p)$ \cite{ErdosRenyi}. We refer to this model as the complete graph angular synchronization model.

The elements of $H$ in the complete graph angular synchronization model are random variables given by the following mixture model.
With probability $p$ the edge $\{i,j\}$ is good and $H_{ij}=e^{\imath (\theta_i - \theta_j)}$, whereas with probability $1-p$ the edge is bad and
$H_{ij}\sim Uniform\left(S^1\right)$. It is convenient to define the diagonal elements as $H_{ii}=p$.

The matrix $H$ is Hermitian and the expected value of its elements is
\begin{equation}
\mathbb{E}H_{ij} = p\,e^{\imath (\theta_i-\theta_j)}.
\end{equation}
In other words, the expected value of $H$ is the rank-one matrix
\begin{equation}
\mathbb{E}H = np zz^*,
\end{equation}
where $z$ is the normalized vector ($\|z\|=1$) given by
\begin{equation}
\label{z}
z_i = \frac{1}{\sqrt{n}}\,e^{\imath \theta_i},\quad i=1,\ldots,n.
\end{equation}
The matrix $H$ can be decomposed as
\begin{equation}
\label{HR}
H = npzz^* + R,
\end{equation}
where $R=H-\mathbb{E}H$ is a random matrix whose elements have zero mean, with $R_{ii}=0$, and for $i\neq j$
\begin{equation}
R_{ij} = \left\{\begin{array}{cc} (1-p)e^{\imath (\theta_i-\theta_j)} & \mbox{with probability } p \\
e^{\imath\varphi} - pe^{\imath (\theta_i-\theta_j)} & \mbox{w.p. } 1-p \mbox{ and } \varphi \sim Uniform([0,2\pi))
\end{array} \right..
\end{equation}
The variance of $R_{ij}$ is
\begin{equation}
\label{var-R}
\mathbb{E}|R_{ij}|^2 = (1-p)^2 p + (1+p^2)(1-p) = 1-p^2
\end{equation}
for $i\neq j$ and 0 for the diagonal elements. Note that for $p=1$ the variance vanishes as all edges become good.

The distribution of the eigenvalues of the random matrix $R$ follows Wigner's semi-circle law \cite{Wigner1,Wigner2} whose support is
$[-2\sqrt{n(1-p^2)},2\sqrt{n(1-p^2)}]$.
The largest eigenvalue
of $R$, denoted $\lambda_1(R)$, is concentrated near the right edge of the support \cite{AlonVu} and the universality of the edge of the spectrum \cite{Soshnikov} implies that it follows the Tracy-Widom distribution \cite{TracyWidom} even when the entries of $R$ are non-Gaussian. For our purposes, the approximation
\begin{equation}
\label{R}
\lambda_1(R) \approx 2\sqrt{n(1-p^2)}
\end{equation}
will suffice, with the probabilistic error bound given in \cite{AlonVu}.

The matrix $H=npzz^* + R$ can be considered as a rank-one perturbation to a random matrix. The distribution of the largest eigenvalue of such perturbed random matrices was investigated in
\cite{Peche,PecheFeral,Furedi} for the particular case where $z$ is proportional to the all-ones vector $(1 \; 1 \cdots \; 1)^T$. Although our vector $z$ given by (\ref{z}) is different, without loss of generality, we can assume $\theta_1=\theta_2=\ldots=\theta_n = 0$, since this assumption does not change the statistical properties of the random matrix $R$. Thus, adopting \cite[Theorem 1.1]{PecheFeral} to $H$ gives that for
\begin{equation}
\label{cond-HR}
np > \sqrt{n(1-p^2)}
\end{equation}
the largest eigenvalue $\lambda_1(H)$ jumps outside the support of the semi-circle law and is normally distributed with mean $\mu$ and variance $\sigma^2$ given by
\begin{equation}
\label{lambda-H}
\lambda_1(H) \sim \mathcal{N}(\mu,\sigma^2),\quad \mu=\frac{np}{\sqrt{1-p^2}} + \frac{\sqrt{1-p^2}}{p}, \quad \sigma^2 = \frac{(n+1)p^2-1}{np^2}(1-p^2),
\end{equation}
whereas for $np<\sqrt{n(1-p^2)}$, $\lambda_1(H)$ still tends to the right edge of the semicircle given at
$2\sqrt{n(1-p^2)}$.

Note that the factor of 2 that appears in (\ref{R}) has disappeared from (\ref{cond-HR}), which is perhaps somewhat non-intuitive: it is expected that $\lambda_1(H) > \lambda_1(R)$ whenever $np > \lambda_1(R)$, but the theorem guarantees that $\lambda_1(H)>\lambda_1(R)$ already for $ np > \frac{1}{2} \lambda_1(R)$.

The condition (\ref{cond-HR}) also implies a lower bound on the correlation between the normalized top eigenvector $v_1$ of $H$ and the vector $z$. To that end, consider the eigenvector equation satisfied by $v_1$:
\begin{equation}
\label{eig-eq}
\lambda_1(H) v_1 = Hv_1 = (np zz^* + R)v_1.
\end{equation}
Taking the dot product with $v_1$ yields
\begin{equation}
\lambda_1(H) = np\,|\langle z, v_1 \rangle|^2 + v_1^* R v_1.
\end{equation}
From $v_1^* R v_1 \leq \lambda_1(R)$ we obtain the lower bound
\begin{equation}
|\langle z, v_1 \rangle|^2 \geq \frac{\lambda_1(H) - \lambda_1(R)}{np},
\end{equation}
with $\lambda_1(H)$ and $\lambda_1(R)$ given by (\ref{R}) and (\ref{lambda-H}). Thus, if the spectral gap $\lambda_1(H)-\lambda_1(R)$ is large enough then $v_1$ must be close to $z$, in which case the eigenvector method successfully recovers the unknown angles. Since the variance of the correlation of two random unit vectors in $\mathbb{R}^n$ is $1/n$, the eigenvector method would give above random correlation values whenever
\begin{equation}
\label{cond-10}
\frac{\lambda_1(H) - \lambda_1(R)}{np} > \frac{1}{n}.
\end{equation}
Replacing in (\ref{cond-10}) $\lambda_1(H)$ by $\mu$ from (\ref{lambda-H}) and $\lambda_1(R)$ by (\ref{R}) and multiplying by $p\sqrt{n}$ yields the condition
\begin{equation}
\label{cond-11}
\frac{\sqrt{n}p}{\sqrt{1-p^2}} + \frac{\sqrt{1-p^2}}{\sqrt{n}p} - 2\sqrt{1-p^2} > \frac{p}{\sqrt{n}}.
\end{equation}
Since $\frac{\sqrt{n}p}{\sqrt{1-p^2}} + \frac{\sqrt{1-p^2}}{\sqrt{n}p} \geq 2$, it follows that (\ref{cond-11}) is satisfied for
\begin{equation}
\label{pn}
p > \frac{1}{\sqrt{n}}.
\end{equation}
Thus, already for $p>\frac{1}{\sqrt{n}}$ we should obtain above random correlations between the vector of angles $z$ and the top eigenvector $v_1$. We therefore define the threshold probability $p_c$ as
\begin{equation}
\label{pcn}
p_c = \frac{1}{\sqrt{n}}.
\end{equation}

When $np \gg \lambda_1(R)$, the correlation between $v_1$ and $z$ can be predicted by using regular perturbation theory for solving the eigenvector equation (\ref{eig-eq}) in an asymptotic expansion with the small parameter $\epsilon = \frac{\lambda_1(R)}{np}$. Such perturbations are derived in standard textbooks on quantum mechanics aiming to find approximations to the energy levels and eigenstates of perturbed time-independent Hamiltonians (see, e.g., \cite[Chapter 6]{Griffiths}). In our case, the resulting asymptotic expansions of the non-normalized eigenvector $v_1$ and of the eigenvalue $\lambda_1(H)$ are given by
\begin{equation}
\label{v1}
v_1 \sim z + \frac{1}{np}\left[Rz - (z^*Rz)z \right] + \ldots,
\end{equation}
and
\begin{equation}
\lambda_1(H) \sim np + z^* R z + \ldots.
\end{equation}
Note that the first order term in (\ref{v1}) is perpendicular to the leading order term $z$, from which it follows that the angle $\alpha$ between the eigenvector $v_1$ and the vector of true angles $z$ satisfies the asymptotic relation
\begin{equation}
\label{tan-alpha}
\tan ^2 \alpha \sim \frac{\|Rz\|^2 - (z^* R z)^2}{(np)^2} + \ldots,
\end{equation}
because $\|Rz - (z^*Rz)z \|^2 = \| Rz \|^2 - (z^* Rz)^2$. The expected values of the numerator terms in (\ref{tan-alpha}) are given by
\begin{equation}
\mathbb{E} \|Rz\|^2 = \mathbb{E} \sum_{i=1}^n \left|\sum_{j=1}^n R_{ij}z_j \right|^2 = \sum_{i,j=1}^n \operatorname{Var}(R_{ij} z_j) = \sum_{i=1}^n \sum_{j\neq i} |z_j|^2 (1-p^2) = (n-1)(1-p^2), \label{term1}
\end{equation}
and
\begin{eqnarray}
\mathbb{E} (z^* R z)^2 &=& \mathbb{E} \left[\sum_{i,j=1}^n R_{ij} \bar{z}_i z_j \right]^2 = \sum_{i,j=1}^n \operatorname{Var} (R_{ij} \bar{z}_i z_j) = (1-p^2)\sum_{i\neq j} |z_i|^2|z_j|^2 \nonumber \\
&=& (1-p^2)\left[\left(\sum_{i=1}^n |z_i|^2 \right)^2 - \sum_{i=1}^n |z_i|^4\right] = (1-p^2)\left(1-\frac{1}{n} \right),\label{term2}
\end{eqnarray}
where we used that $R_{ij}$ are i.i.d zero mean random variables with variance given by (\ref{var-R}) and that $|z_i|^2 = \frac{1}{n}$. Substituting (\ref{term1})-(\ref{term2}) into (\ref{tan-alpha}) results in
\begin{equation}
\mathbb{E} \tan^2 \alpha \sim \frac{(n-1)^2 (1-p^2)}{n^3 p^2} + \ldots,
\end{equation}
which for $p \ll 1$ and $n \gg 1$ reads
\begin{equation}
\mathbb{E} \tan^2 \alpha \sim \frac{1}{np^2} + \ldots.
\end{equation}
This expression shows that as $np^2$ goes to infinity, the angle between $v_1$ and $z$ goes to zero and the correlation between them goes to 1. For $np^2 \gg 1$, the leading order term in the expected squared correlation $\mathbb{E}\cos^2 \alpha$ is given by
\begin{equation}
\label{cos-alpha}
\mathbb{E}\cos^2 \alpha = \mathbb{E}\frac{1}{1+\tan^2 \alpha} \sim \frac{1}{1 +\displaystyle\frac{1}{np^2}} + \ldots.
\end{equation}
We conclude that even for very small $p$ values, the eigenvector method successfully recovers the angles if there are enough equations, that is, if $np^2$ is large enough.

\begin{figure}[h]
\begin{centering}
\subfigure[$p=0.15$]{
\includegraphics[width=0.3\textwidth]{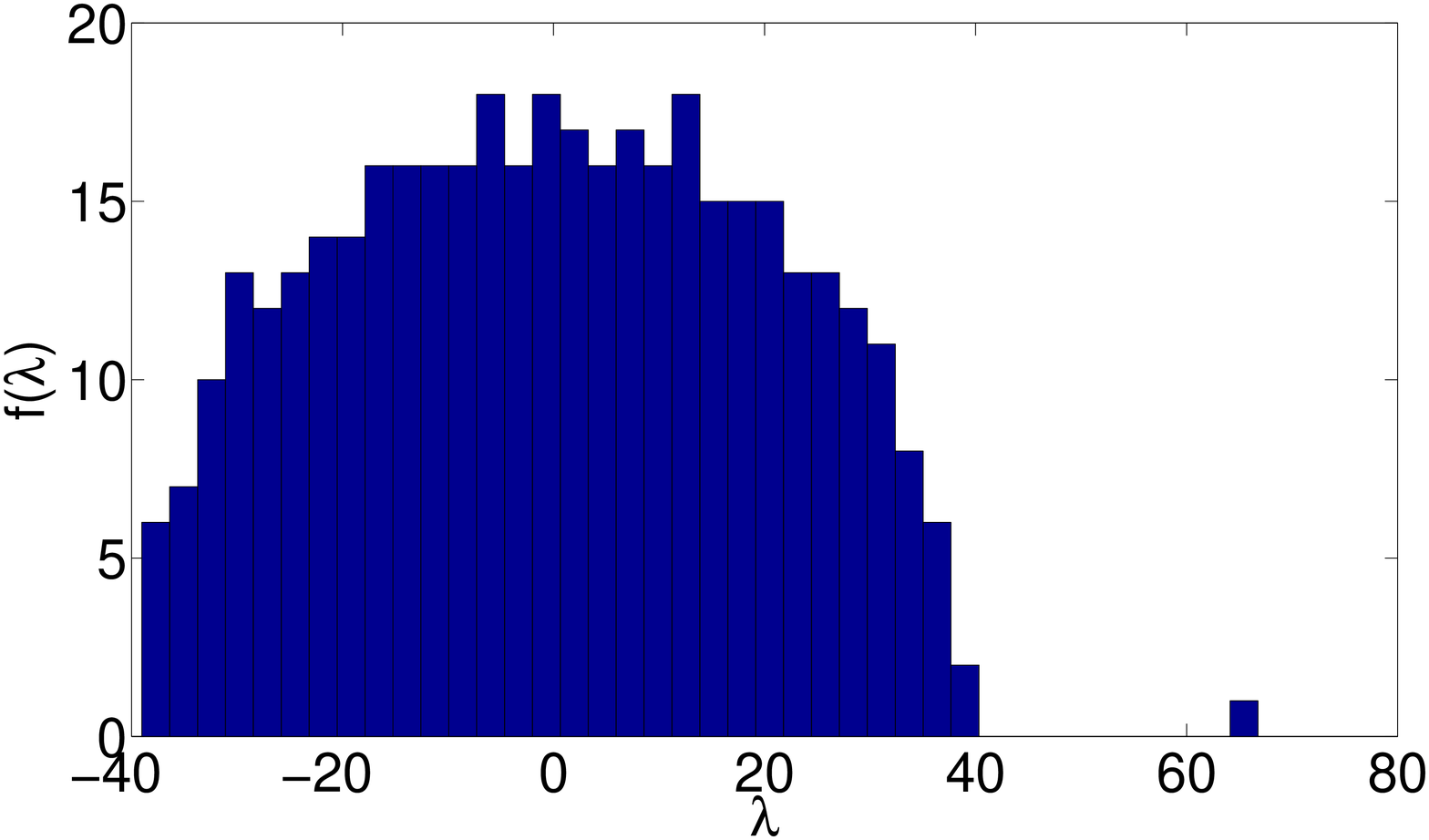}
\label{fig:N400p15}
}
\subfigure[$p=0.1$]{
\includegraphics[width=0.3\textwidth]{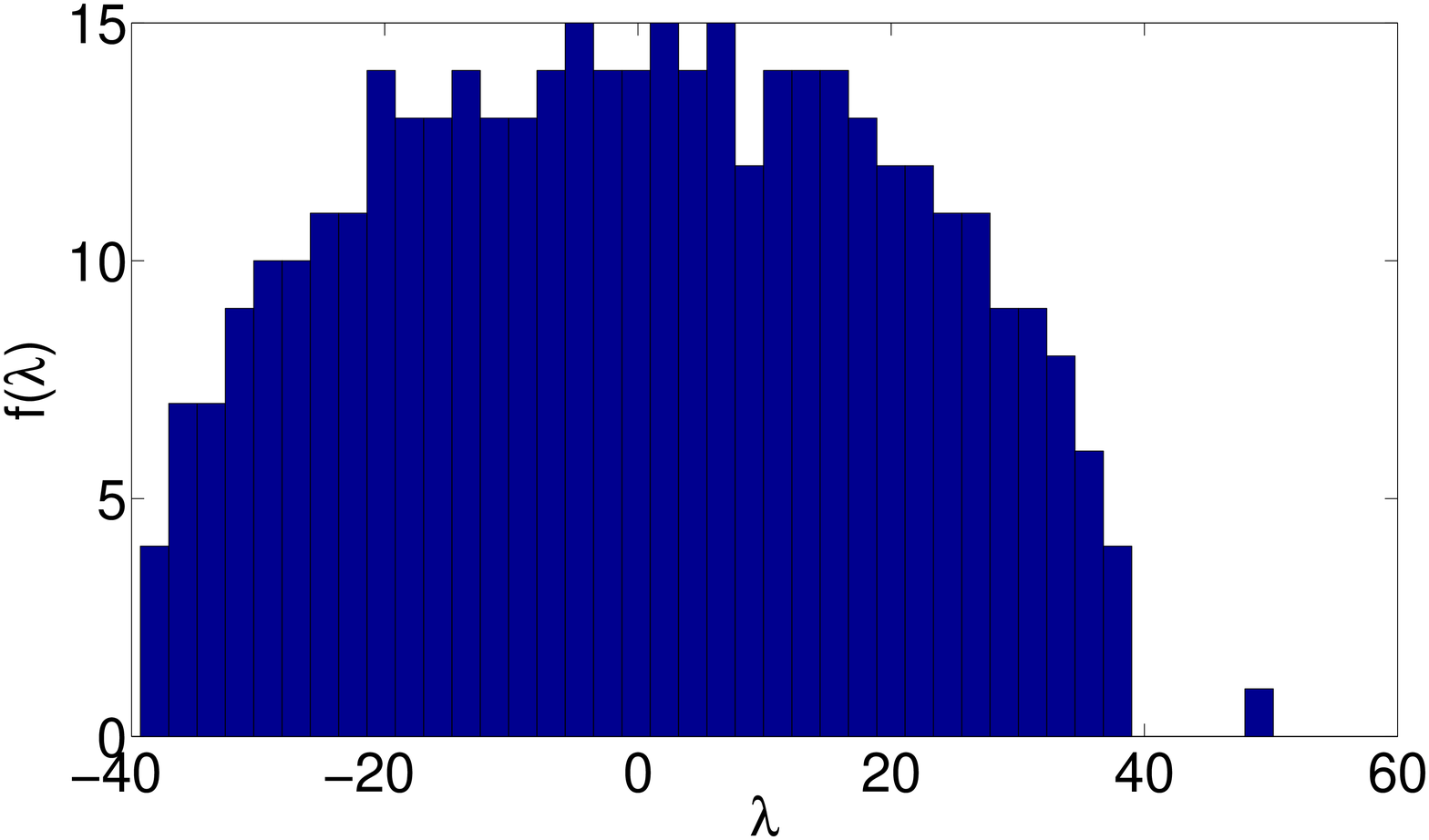}
\label{fig:N400p10}
}
\subfigure[$p=0.05$]{
\includegraphics[width=0.3\textwidth]{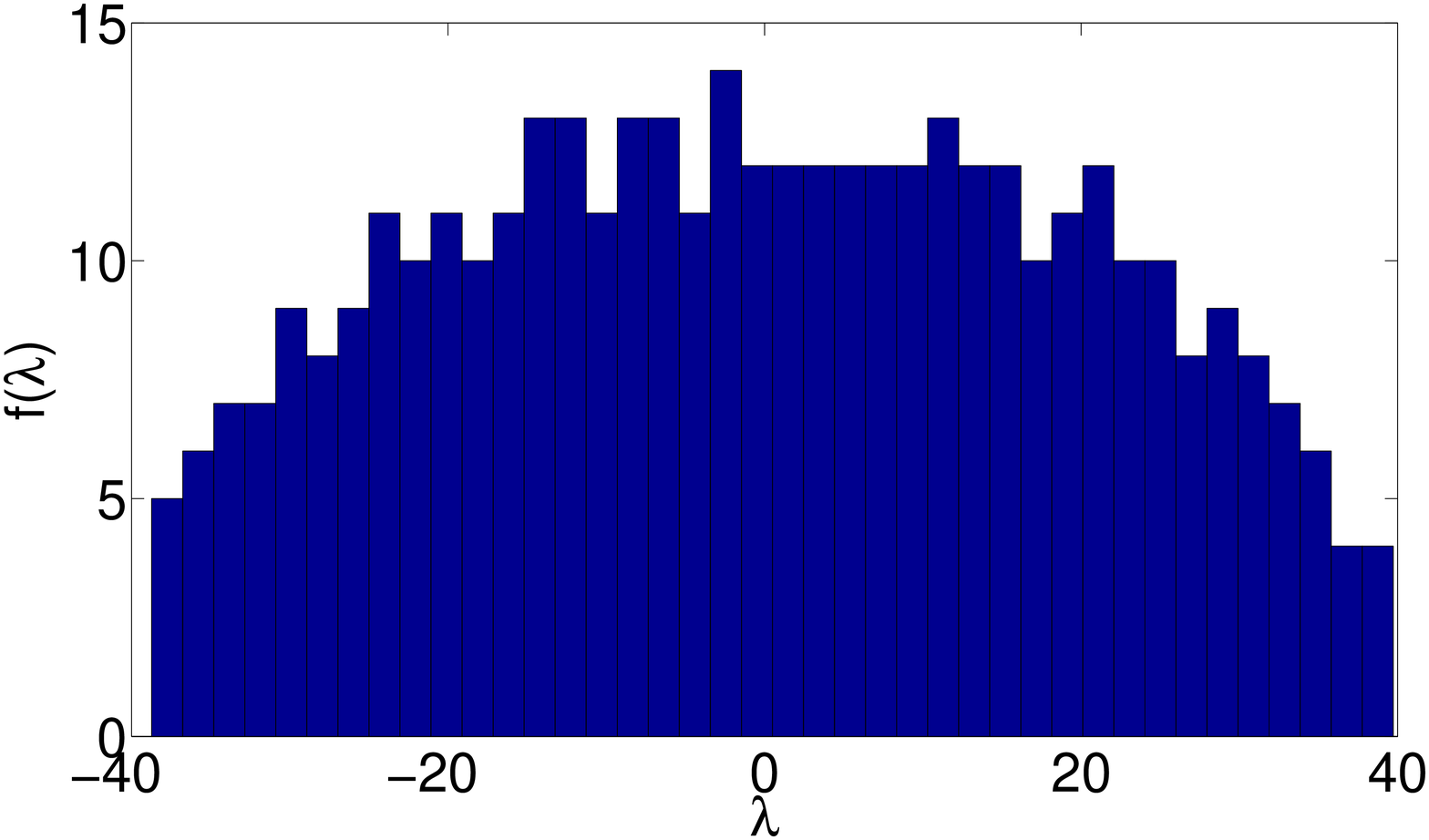}
\label{fig:N400p05}
}
\end{centering}
\caption{Histogram of the eigenvalues of the matrix $H$ in the complete graph model for $n=400$ and different values of $p$.}\label{fig:N400}
\end{figure}

Figure \ref{fig:N400} shows the distribution of the eigenvalues of the matrix $H$ for $n=400$ and different values of $p$.
The spectral gap decreases as $p$ is getting smaller. From (\ref{cond-HR}) we expect a spectral gap for $p\geq p_c$ where the critical value is $p_c=\frac{1}{\sqrt{400}}=0.05$. The experimental values of $\lambda_1(H)$ also agree with (\ref{lambda-H}). For example, for $n=400$ and $p=0.15$, the expected value of the largest eigenvalue is $\mu=67.28$ and its standard deviation is
$\sigma = 0.93$, while for $p=0.1$ we get $\mu=50.15$ and
$\sigma = 0.86$; these value are in full agreement with the location of the largest eigenvalues in Figures \ref{fig:N400p15}-\ref{fig:N400p10}. Note that the right edge of the semi-circle is smaller than $2\sqrt{n} = 40$, so the spectral gap is significant even when $p=0.1$.

The skeptical reader may wonder whether the existence of a visible spectral gap necessarily implies that the normalized top eigenvector $v_1$
correctly recovers the original set of angles $\theta_1,\ldots,\theta_n$ (up to a constant phase). To that end, we compute the following two measures of correlation $\rho_1$ and $\rho_2$ for the correlation between the vector of true angles $z$ and the computed normalized top eigenvector $v_1$:
\begin{equation}
\rho_1 = \left|\frac{1}{n}\sum_{i=1}^n e^{-i\theta_i} \frac{v_1(i)}{|v_1(i)|}\right|,\quad \rho_2 = \left|\frac{1}{\sqrt{n}}\sum_{i=1}^n e^{-i\theta_i} v_1(i) \right| = |\langle z,v_1 \rangle |.
\end{equation}
The correlation $\rho_1$ takes into account the rounding procedure (\ref{eig2}), while $\rho_2$ is simply the dot product between $v_1$ and $z$ without applying any rounding. Clearly, $\rho_1,\rho_2 \leq 1$ (Cauchy-Schwartz), and $\rho_1=1$ iff the two sets of angles are the same up to a rotation. Note that it is possible to have $\rho_1 = 1$ with $\rho_2 < 1$. This happens when the angles implied by $v_1(i)$ are all correct, but the magnitudes $|v_1(i)|$ are not all the same. Table \ref{tab:eig-comp} summarizes the experimentally obtained correlations $\rho_1,\rho_2$ for different values of $p$ with $n=100$ (Table \ref{tab:t1}) and $n=400$ (Table \ref{tab:t2}). The experimental results show that for large values of $np^2$ the correlation is very close to 1, indicating a successful recovery of the angles. The third column, indicating the values of $\left(1+\frac{1}{np^2}\right)^{-\frac{1}{2}}$ is  motivated by the asymptotic expansion (\ref{cos-alpha}) and seems to provide a very good approximation for $\rho_2$ when $np^2 \gg 1$, with deviations attributed to higher order terms of the asymptotic expansion and to statistical fluctuations around the mean value. Below the threshold probability (ending rows of Tables \ref{tab:t1} and \ref{tab:t2} with $np^2 < 1$), the correlations take values near $\frac{1}{\sqrt{n}}$, as expected from the correlation of two unit random vectors in $\mathbb{R}^n$ ($\frac{1}{\sqrt{100}} = 0.1$ and $\frac{1}{\sqrt{400}} = 0.05$).

From the practical point of view, most important is the fact that the eigenvector method successfully recovers the angles even when a large portion of the offset measurements consists of just outliers. For example, for $n=400$, the correlation obtained when $85\%$ of the offset measurements were outliers (only $15\%$ are good measurements) was $\rho_1=0.97$.

\begin{table}[h]
\begin{centering}
\subfigure[$n=100$]{
\label{tab:t1}
\begin{tabular}{|c|c|c|c|c|}
\hline
$p$ & $np^2$ & $\left(1+\frac{1}{np^2}\right)^{-\frac{1}{2}}$ & $\rho_1$ & $\rho_2$ \\
\hline
0.4 & 16 & 0.97 & 0.99 & 0.98 \\
0.3 & 9 & 0.95 & 0.97 & 0.95 \\
0.2 & 4 & 0.89 & 0.90 & 0.88 \\
0.15 & 2.25 & 0.83 & 0.75 & 0.81 \\
0.1 & 1 & 0.71 &  0.34 & 0.35 \\
0.05 & 0.25 & 0.45 &  0.13 & 0.12\\
\hline
\end{tabular}}
\subfigure[$n=400$]{
\label{tab:t2}
\begin{tabular}{|c|c|c|c|c|}
\hline
$p$ & $np^2$ & $\left(1+\frac{1}{np^2}\right)^{-\frac{1}{2}}$ & $\rho_1$ & $\rho_2$ \\
\hline
0.2 & 16 & 0.97 & 0.99 & 0.97 \\
0.15& 9 & 0.95 & 0.97  & 0.95 \\
0.1 & 4 & 0.89 & 0.90 &  0.87 \\
0.075 & 2.25 & 0.83 & 0.77 & 0.76\\
0.05 & 1 & 0.71 & 0.28 & 0.32 \\
0.025 & 0.25 & 0.45 & 0.06 & 0.07 \\
\hline
\end{tabular}}
\end{centering}
\caption{Correlations between the top eigenvector $v_1$ of $H$ and the vector $z$ of true angles for different values of $p$ in the complete graph model.} \label{tab:eig-comp}
\end{table}

\subsection{Analysis of the angular synchronization problem in general}
We turn to analyze the eigenvector method for general measurement graphs, where the graph of good measurements is assumed to be connected, while the graph of bad edges is assumed to be made of edges that are uniformly drawn from the remaining edges of the complete graph once the good edges has been removed from it.
Our analysis is based on generalizing the decomposition given in (\ref{HR}).

Let $A$ be the adjacency matrix for the set of good edges $E_{good}$:
\begin{equation}
A_{ij} = \left\{\begin{array}{cc} 1 & \{i,j\} \in E_{good} \\
0 & \{i,j\} \not \in E_{good}
\end{array} \right..
\end{equation}
As the matrix $A$ is symmetric, it has a complete set of real eigenvalues $\lambda_1 \geq \lambda_2 \geq \ldots \geq \lambda_n$ and corresponding real orthonormal eigenvectors
$\psi_1,\ldots,\psi_n$ such that
\begin{equation}
A = \sum_{l=1}^n \lambda_l \psi_l \psi_l^T.
\end{equation}

Let $Z$ be an $n\times n$ diagonal matrix whose diagonal elements are $Z_{ii} = e^{\imath \theta_i}$. Clearly, $Z$ is a unitary matrix ($ZZ^*=I$). Define the Hermitian matrix $B$ by conjugating $A$ with $Z$
\begin{equation}
\label{B}
B = Z A Z^*.
\end{equation}
It follows that the eigenvalues of $B$ are equal to the eigenvalues $\lambda_1,\ldots,\lambda_n$ of $A$, and the corresponding eigenvectors $\{\phi_l\}_{l=1}^n$ of $B$, satisfying $B\phi_l = \lambda_l \phi_l$ are given by
\begin{equation}
\phi_l = Z \psi_l, \quad l=1,\ldots,n.
\end{equation}
From (\ref{B}) it follows that
\begin{equation}
B_{ij} = \left\{\begin{array}{cc} e^{\imath (\theta_i - \theta_j)} & \{i,j\}\in E_{good} \\
0 & \{i,j\}\not \in E_{good}
\end{array} \right..
\end{equation}
We are now ready to decompose the matrix $H$ defined in (\ref{H}) as
\begin{equation}
\label{HBR}
H = B + R,
\end{equation}
where $R$ is a random matrix whose elements are given by
\begin{equation}
\label{Rsparse}
R_{ij} = \left\{\begin{array}{cc} e^{\imath \delta_{ij}} & \{i,j\}\in E_{bad} \\
0 & \{i,j\}\not \in E_{bad}
\end{array} \right.,
\end{equation}
where $\delta_{ij}\sim Uniform([0,2\pi)$ for $\{i,j\}\in E_{bad}$.
The decomposition (\ref{HBR}) is extremely useful, because it sheds light into the eigen-structure of $H$ in terms of the much simpler eigen-structures of $B$ and $R$.

First, consider the matrix $B$ defined in (\ref{B}), which shares the same spectrum with $A$ and whose eigenvectors $\phi_1,\ldots,\phi_n$
are phase modulations of the eigenvectors $\psi_1,\ldots,\psi_n$ of $A$. If the graph of good measurements is connected, as it must be in order to have a unique solution for the angular synchronization problem (see second paragraph of Section \ref{sec:intro}), then the Perron-Frobenius theorem (see, e.g., \cite[Chapter 8]{HornMatrix}) for the non-negative matrix $A$ implies that the entries of $\psi_1$ are all positive
\begin{equation}
\label{positive}
\psi_1(i) > 0, \quad \mbox{for all } i=1,2,\ldots,n,
\end{equation}
and therefore the complex phases of the coordinates of the top eigenvector $\phi_1=Z\psi_1$ of $B$ are identical to the true angles, that is, $e^{\imath \theta_i} = \frac{\phi_1(i)}{|\phi_1(i)|}$. Hence, if the top eigenvector of $H$ is highly correlated with the top eigenvector of $B$, then the angles will be estimated with high accuracy. We will shortly derive the precise condition that guarantees such a high correlation between the eigenvectors of $H$ and $B$.

The spectral gap $\Delta_{good}$ of the good graph is the difference between its first and second eigenvalues, i.e., $\Delta_{good} = \lambda_1(A)-\lambda_2(A)$. The Perron-Frobenius theorem and the connectivity of the graph of good measurements also imply that $\Delta_{good} > 0$.

Next we turn to analyze the spectrum of the random matrix $R$ given in (\ref{Rsparse}). We assume that the $m_{bad}$ bad edges were drawn uniformly at random from the remaining edges of the complete graph on $n$ vertices that are not already good edges. There are only $2m_{bad}$ nonzero elements in $R$, which makes $R$ a sparse matrix with an average number of $2m_{bad}/n$ nonzero entries per row. The nonzero entries of $R$ have zero mean and unit variance. The spectral norm of such sparse random matrices was studied in \cite{Khorunzhy2001,Khorunzhiy2003} where it was shown that with probability 1, $$\limsup_{n\to \infty} \frac{\sqrt{n}}{\sqrt{2m_{bad}}}\lambda_1(R) \leq 2$$ as long as $\frac{m_{bad}}{n\log n}\to \infty$ as $n\to \infty$. The implication of this result is that we can approximate $\lambda_1(R)$ with
\begin{equation}
\label{R-bad}
\lambda_1(R) \approx 2\frac{\sqrt{2m_{bad}}}{\sqrt{n}}.
\end{equation}
Similar to the spectral gap condition (\ref{cond-HR}), requiring
\begin{equation}
\label{new-cond}
\Delta_{good} > \frac{1}{2}\lambda_1(R),
\end{equation}
ensures that with high probability, the top eigenvector of $H$ would be highly correlated with the top eigenvector of $B$. Plugging (\ref{R-bad}) into (\ref{new-cond}), we get the condition
\begin{equation}
\label{cond-mmm}
\Delta_{good} > \frac{\sqrt{2m_{bad}}}{\sqrt{n}}.
\end{equation}

We illustrate the above analysis for the small world graph, starting with a neighborhood graph on the unit sphere $S^2$ with $n$ vertices corresponding to points on the sphere and $m$ edges, and rewiring each edge with probability $1-p$ at random, resulting shortcut edges. The shortcut edges are considered as bad edges, while unperturbed edges are the good edges.
As the original $m$ edges of the small world graph are rewired with probability $1-p$, the expected number of bad edges $\mathbb{E}m_{bad}$ and the expected number of good edges $\mathbb{E}m_{good}$ are given by
\begin{equation}
\label{m-exp}
\mathbb{E}m_{good}=pm, \quad \mathbb{E}m_{bad} = (1-p)m,
\end{equation}
with relatively small fluctuations of $O(\sqrt{mp(1-p)})$.

The average degree of the original unperturbed graph is $\bar{d} = \frac{2m}{n}$. Assuming uniform sampling of points on the sphere, it follows that the average area of the spherical cap covered by the neighboring points is $4\pi \frac{\bar{d}}{n} = \frac{8\pi m}{n^2}$. The average opening angle $\eta$ corresponding to this cap satisfies $2\pi (1-\cos \eta) =\frac{8\pi m}{n^2}$, or $1-\cos \eta = \frac{4m}{n^2}$. Consider the limit $m,n\to \infty$ while keeping the ratio $c=4m/n^2$ constant. By the law of large numbers, the matrix $\frac{1}{n}A$ converges in this limit to the integral convolution operator $\mathcal{K}$ on $S^2$ (see, e.g., \cite{Coifman2006a}), given by
\begin{equation}
(\mathcal{K}f)(\beta) = \frac{p}{4\pi }\int_{S^2} \chi_{[1-c, 1]}(\langle \beta,\beta' \rangle) f(\beta')\,dS_{\beta'},\quad \beta\in S^2,
\end{equation}
where $\chi_I$ is the characteristic function of the interval $I$.

The classical Funk-Hecke theorem (see, e.g.,
\cite[p. 195]{Natr2001a}) asserts that the spherical harmonics are
the eigenfunctions of convolution operators over the sphere, and the eigenvalues $\lambda_l$ are given by
$$\lambda_l(\mathcal{K}) = \frac{p}{2} \int_{-1}^1  \chi_{[1-c, 1]} P_l(t)\,dt = \frac{p}{2} \int_{1-c}^1  P_l(t)\,dt,$$
and have multiplicities $2l+1$ ($l=0,1,2,\ldots$), where $P_l$ are the Legendre polynomials ($P_0(t)=1, P_1(t)=t, \ldots$). In particular, $\lambda_0(\mathcal{K}) = \frac{pc}{2}$, $\lambda_1(\mathcal{K}) = \frac{pc}{2}(1-\frac{1}{2}c)$, and the spectral gap of $\mathcal{K}$ is $\Delta(\mathcal{K}) = \frac{pc^2}{4}$. The spectral gap of $A$ is approximately
\begin{equation}
\label{Delta-A}
\Delta_{good} \approx n\Delta (\mathcal{K}) = \frac{npc^2}{4} = \frac{4m^2p}{n^3}.
\end{equation}
Plugging (\ref{m-exp}) and (\ref{Delta-A}) into (\ref{cond-mmm}) yields the condition
\begin{equation}
\frac{4m^2p}{n^3} > \frac{\sqrt{2(1-p)m}}{\sqrt{n}},
\end{equation}
which is satisfied for $p > p_c$, where $p_c$ is the threshold probability
\begin{equation}
\label{p-thresh}
p_c \geq \sqrt{\frac{n^5}{8m^3}}.
\end{equation}
We note that this estimate for the threshold probability is far from being tight and can be improved in principle by taking into account the entire spectrum of the good graph rather than just the spectral gap between the top eigenvalues, but we do not attempt to derive tighter bounds here.

\begin{figure}[h]
\begin{centering}
\subfigure[$p=1$]{
\includegraphics[width=0.22\textwidth]{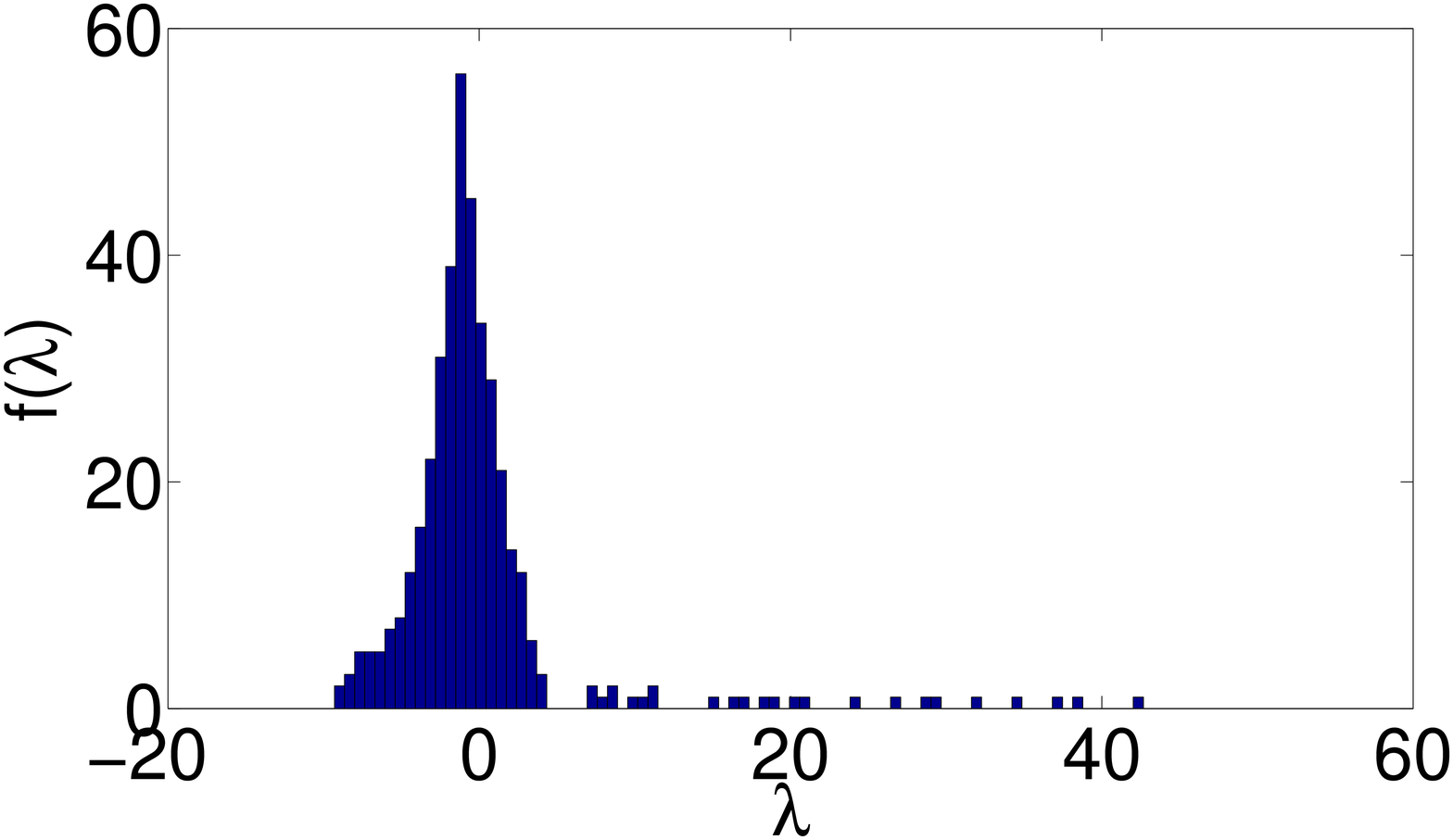}
\label{fig:S2N400p0}
}
\subfigure[$p=0.7$]{
\includegraphics[width=0.22\textwidth]{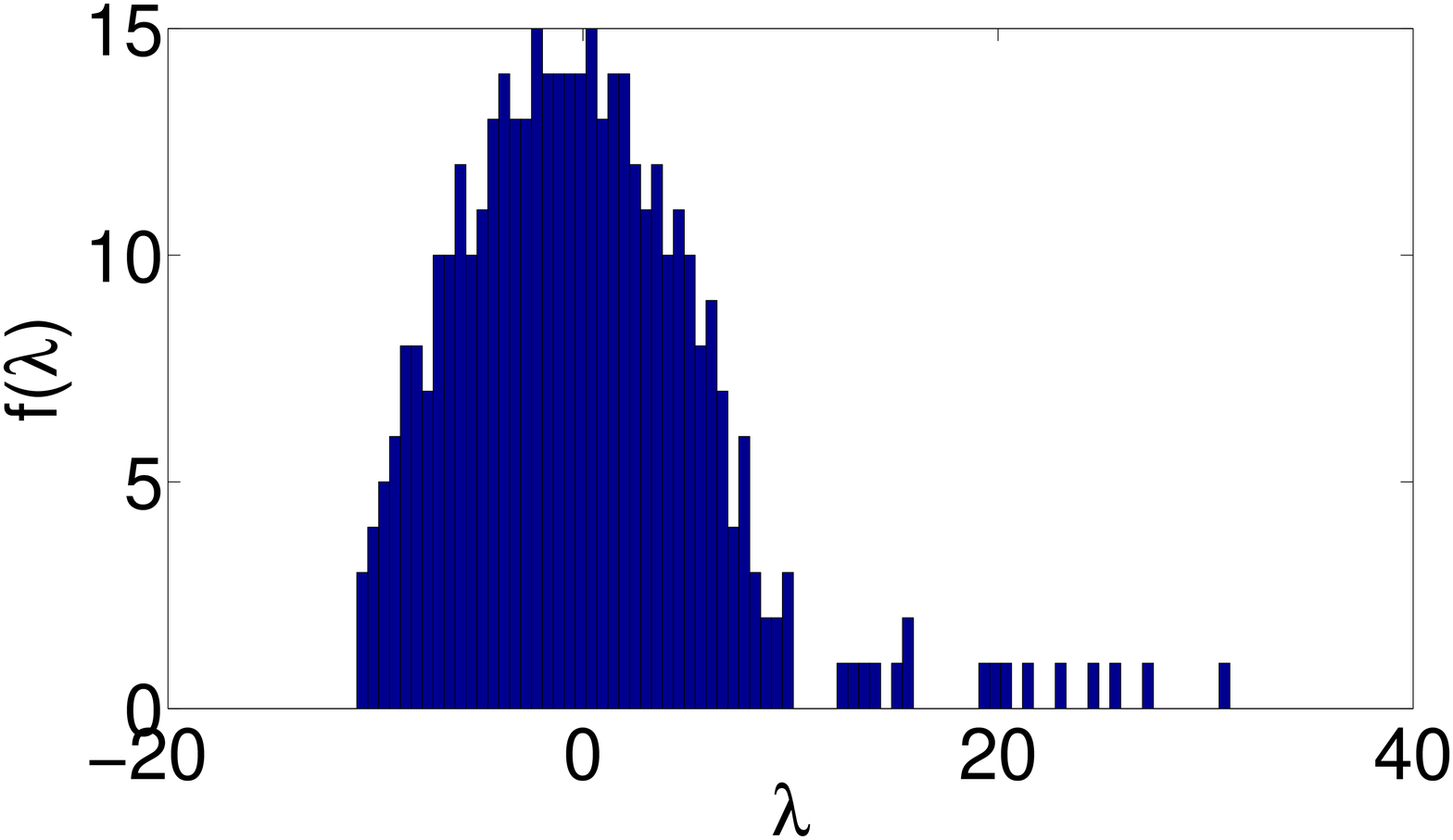}
\label{fig:S2N400p30}
}
\subfigure[$p=0.4$]{
\includegraphics[width=0.22\textwidth]{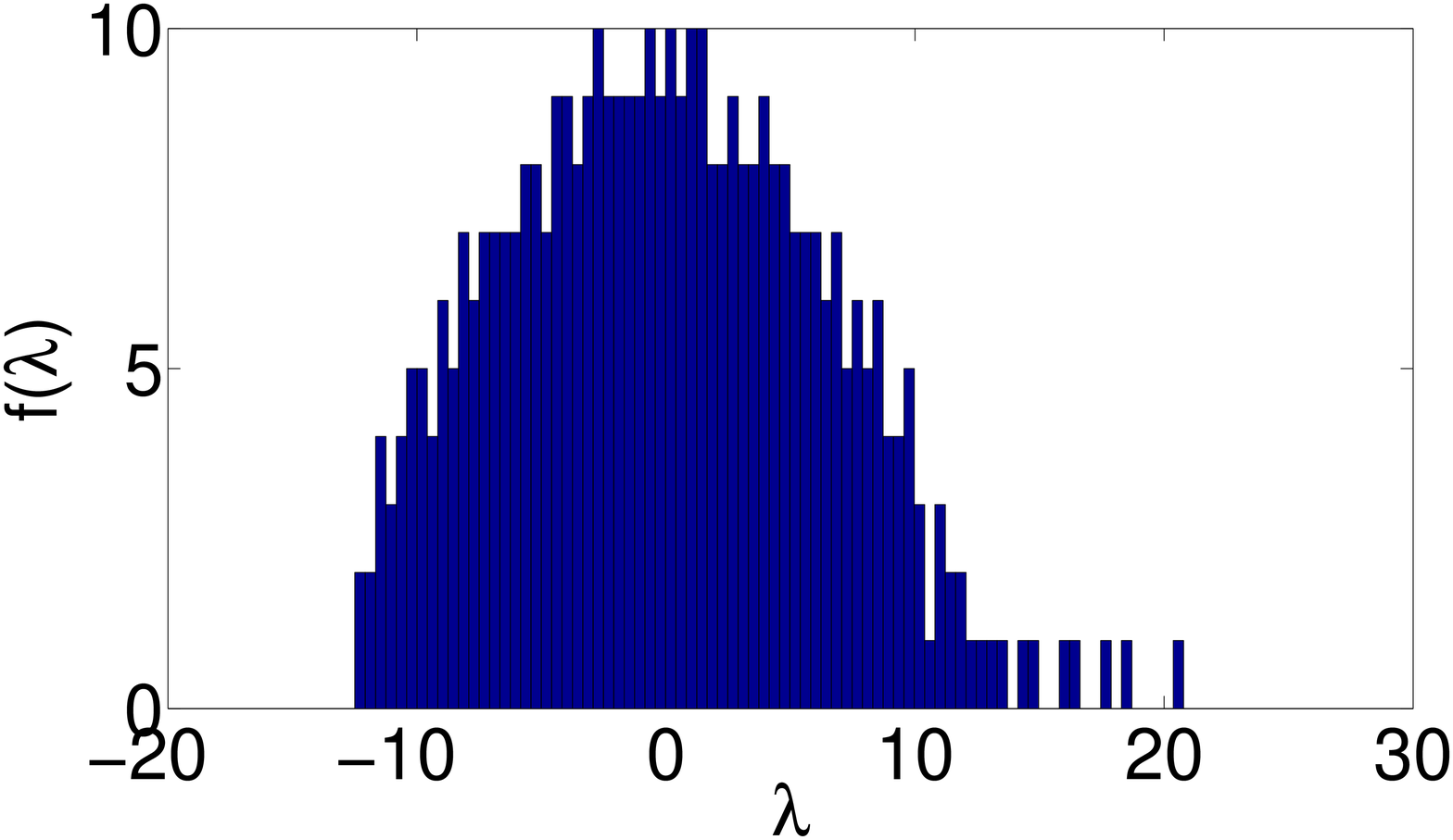}
\label{fig:S2N400p60}
}
\subfigure[$p=0.1$]{
\includegraphics[width=0.22\textwidth]{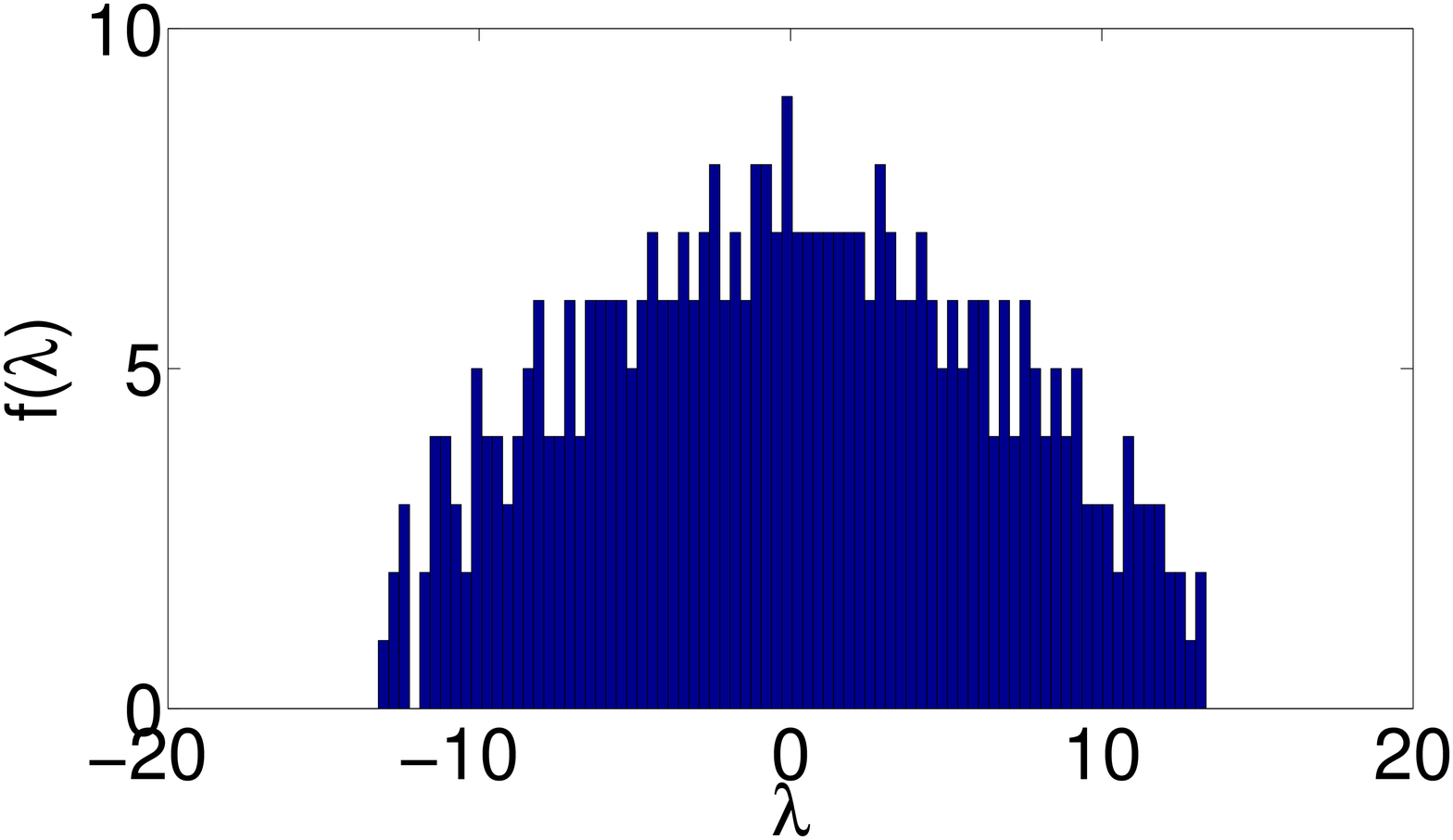}
\label{fig:S2N400p90}
}
\end{centering}
\caption{Histogram of the eigenvalues of the matrix $H$ in the small-world model for $n=400$, $\varepsilon=0.2$, $m\approx 8000$, and different values of $p$.}\label{fig:S2N400}
\end{figure}

We end this section by describing the results of a few numerical experiments. Figure \ref{fig:S2N400} shows the histogram of the eigenvalues of the matrix $H$ for small-world graphs on $S^2$. Each graph was generated by sampling $n$ points $\beta_1,\ldots,\beta_n$
on the unit sphere $S^2$ in $\mathbb{R}^3$ from the uniform distribution as well as $n$ random rotation angles $\theta_1,\ldots,\theta_n$ uniformly distributed in $[0,2\pi)$. An edge between $i$ and $j$ exists iff $\langle \beta_i , \beta_j \rangle > 1-\varepsilon$, where $\varepsilon$ is a small parameter that determines the connectivity (average degree) of the graph. The resulting graph is a neighborhood graph on $S^2$.
The small world graphs were obtained by randomly rewiring the edges of the neighborhood graph. Every edge is rewired with probability $1-p$, so that the expected proportion of good edges is $p$.

\begin{figure}[h]
\begin{centering}
\subfigure[$p=1$]{
\includegraphics[width=0.22\textwidth]{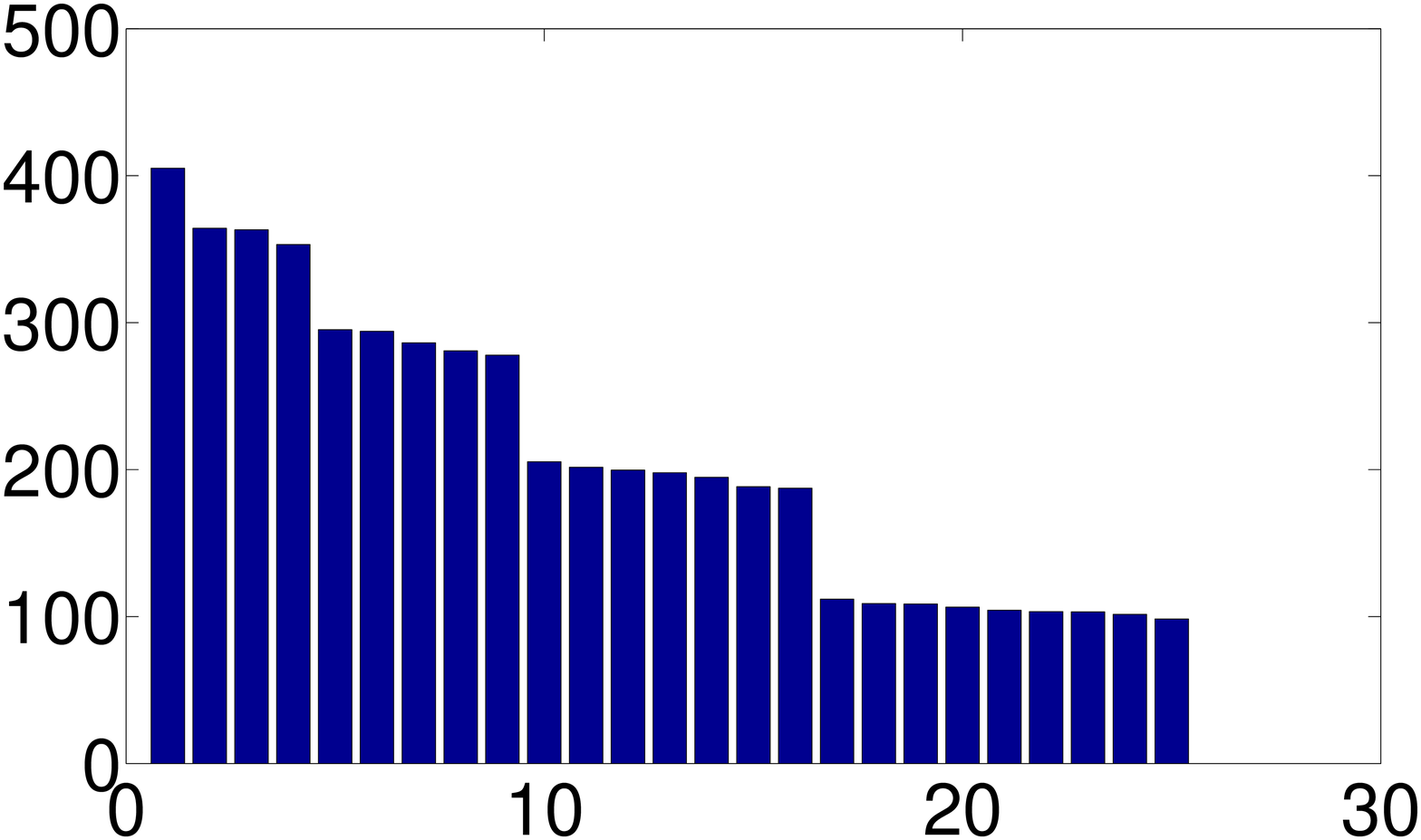}
\label{fig:S2N4000p0}
}
\subfigure[$p=0.2$]{
\includegraphics[width=0.22\textwidth]{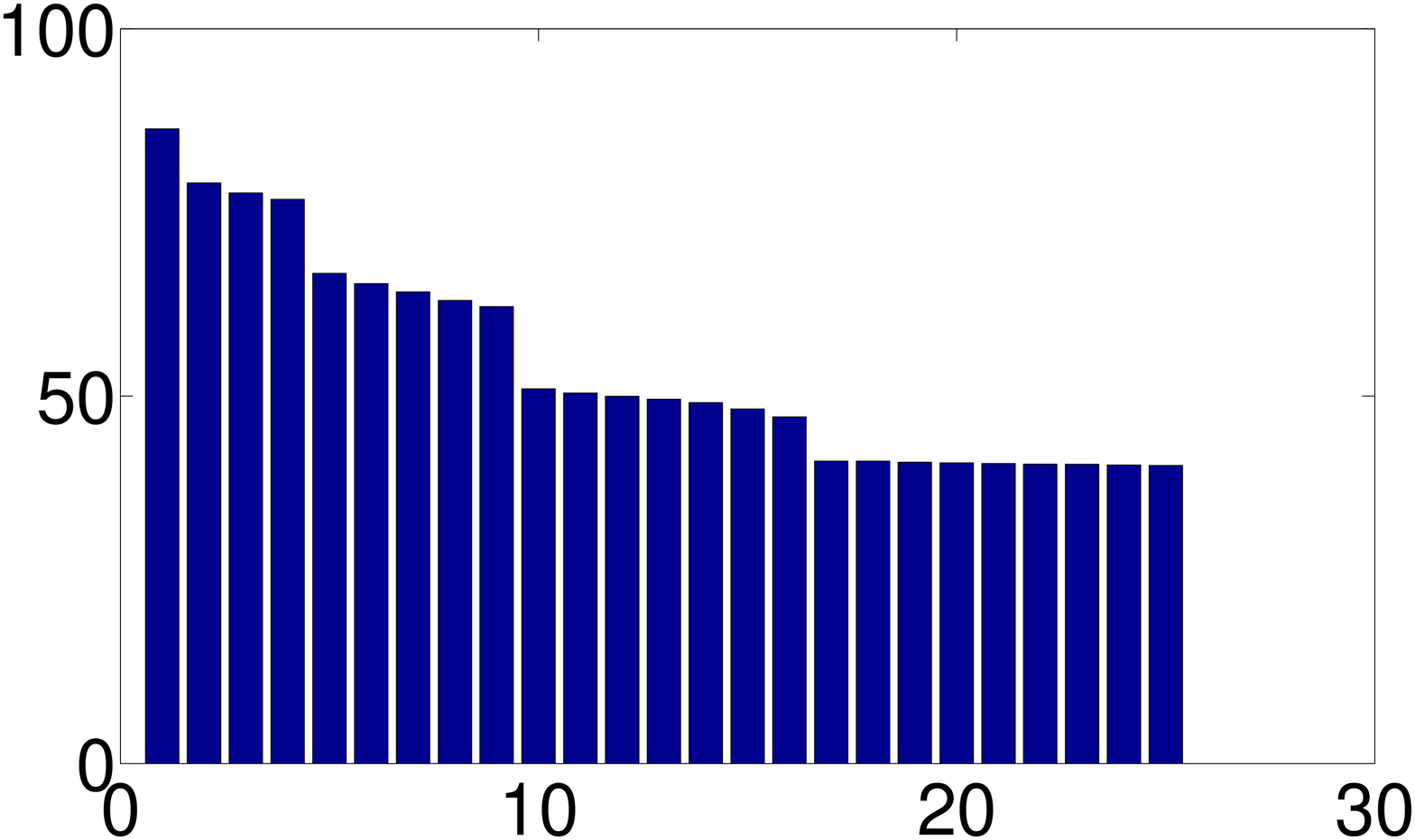}
\label{fig:S2N4000p80}
}
\subfigure[$p=0.1$]{
\includegraphics[width=0.22\textwidth]{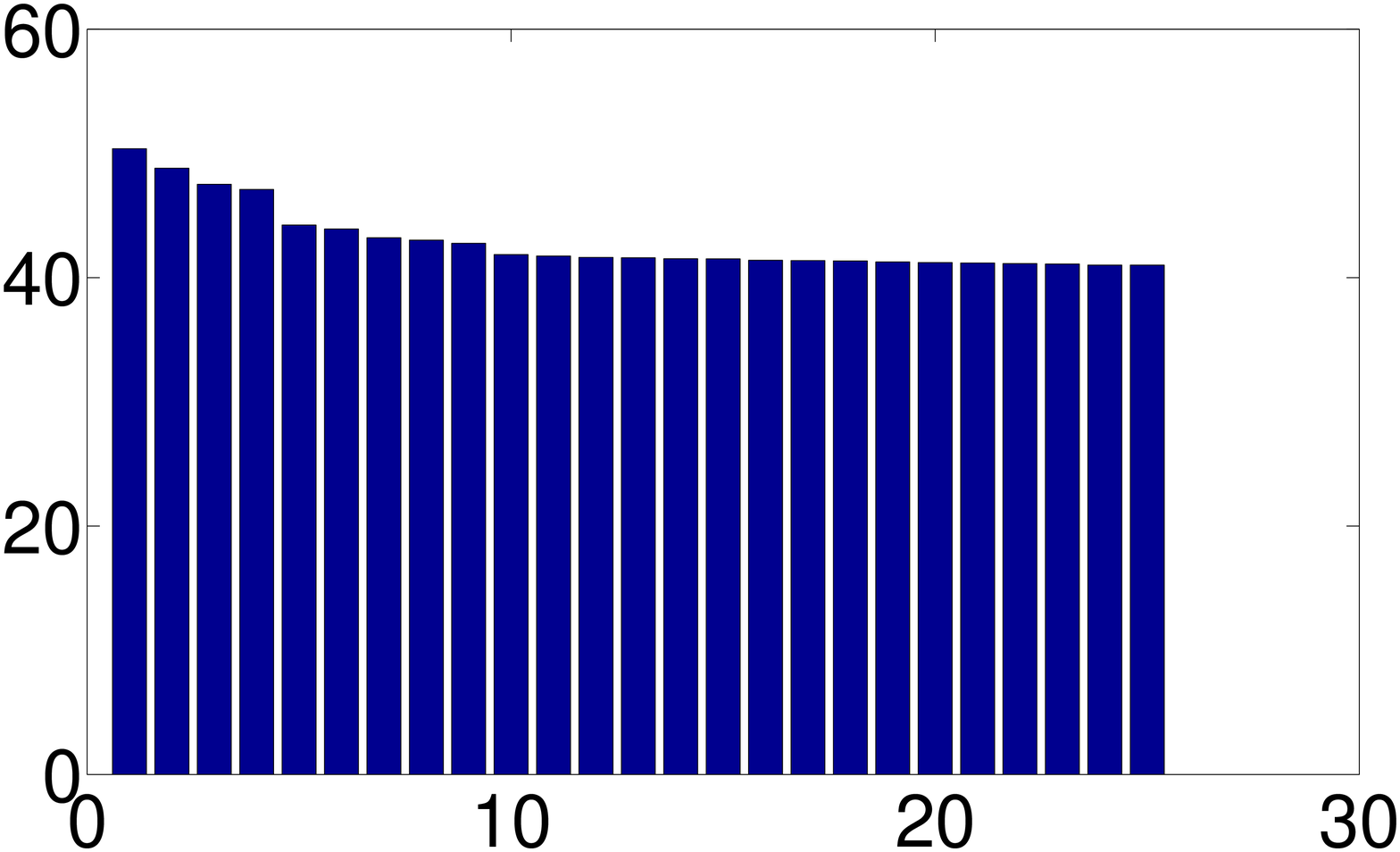}
\label{fig:S2N4000p90}
}
\subfigure[$p=0.05$]{
\includegraphics[width=0.22\textwidth]{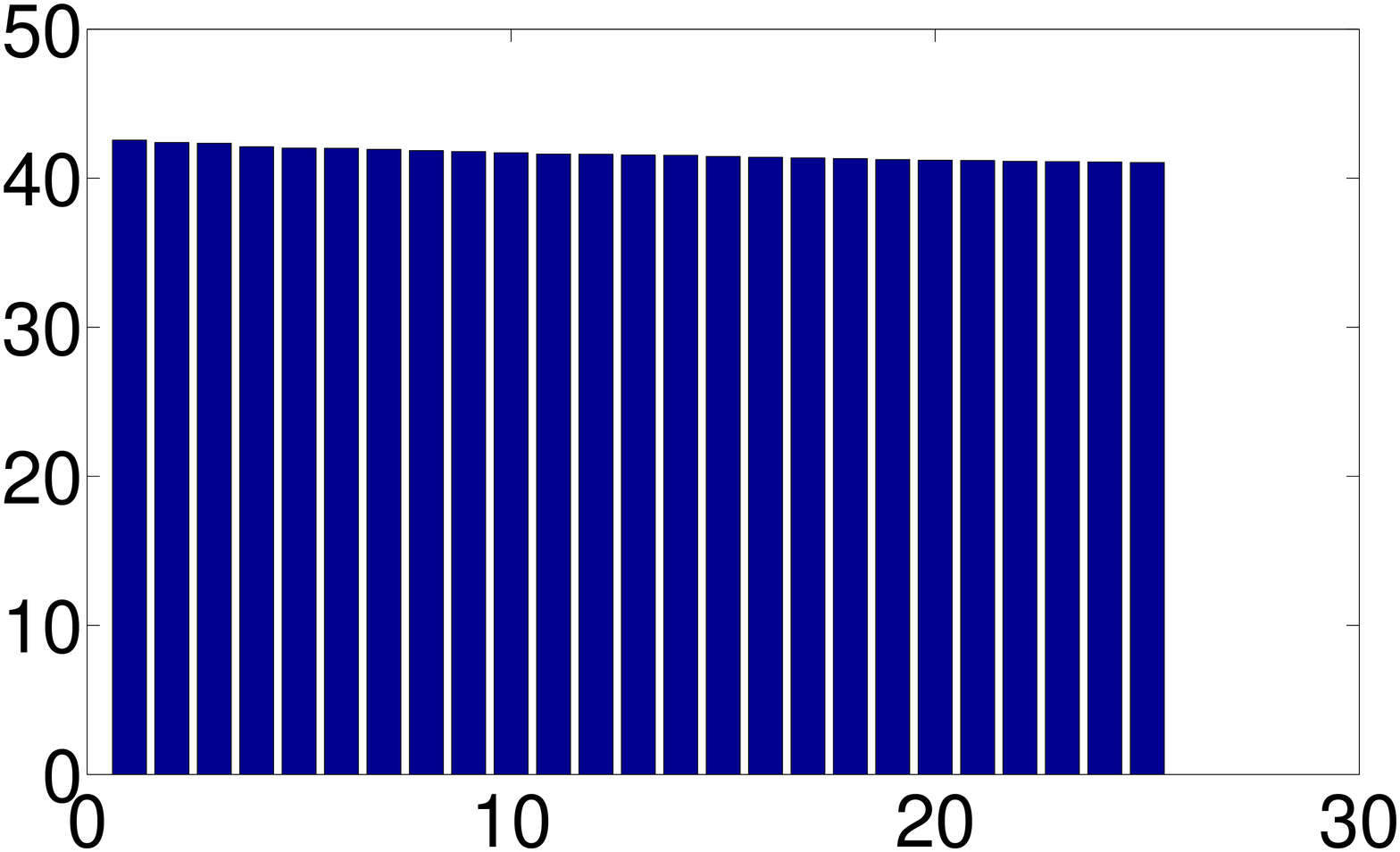}
\label{fig:S2N4000p95}
}
\end{centering}
\caption{Bar plot of the 25 largest eigenvalues of the matrix $H$ in the small-world model for $n=4000$, $\varepsilon=0.2$, $m \approx 8\cdot 10^5$, and different values of $p$. The multiplicities $1,3,5,7,9$ corresponding to the spherical harmonics are evident as long as $p$ is not too small. As $p$ decreases the high-oscillatory spherical harmonics are getting ``swallowed" by the semi-circle.}\label{fig:S2N4000}
\end{figure}

The histograms of Figure \ref{fig:S2N400} for the eigenvalues of $H$ seem to be much more exotic than the ones obtained in the complete graph case shown in Figure \ref{fig:N400}.
In particular, there seems to be a long tail of large eigenvalues, rather than a single eigenvalue that stands out from all the others. But now we understand that these eigenvalues are nothing but the top eigenvalues of the adjacency matrix of the good graph, related to the spherical harmonics. This behavior is better visible in Figure \ref{fig:S2N4000}.

The experimental correlations given in Table \ref{tab:small} indicate jumps in the correlation values that occur between $p=0.15$ and $p=0.2$ for $n=100$ and between $p=0.1$ and $p=0.12$ for $n=400$. The experimental threshold values seem to follow the law $p_c \approx \sqrt{\frac{n}{2m}}$ that holds for the complete graph case (\ref{pn}) with $m={n \choose 2}$. As mentioned earlier, (\ref{p-thresh}) is a rather pessimistic estimate of the threshold probability.

Also evident from Table \ref{tab:small} is that the correlation goes to $1$ as $2mp^2/n \to \infty$. We remark that using regular perturbation theory and the relation of the eigenstructure of $B$ to the spherical harmonics, it should be possible to obtain an asymptotic series for the correlation in terms of the large parameter $2mp^2/n$, similar to the asymptotic expansion (\ref{cos-alpha}).
\begin{table}[h]
\begin{center}
\subfigure[$n=100$, $\varepsilon=0.3$, $m\approx 750$]{
\label{tab:small100}
\begin{tabular}{|c|c|c|}
\hline
$p$ & $\frac{2mp^2}{n}$ & $\rho_1$ \\
\hline
0.8 & 9.6 & 0.923 \\
0.6 & 5.4 & 0.775 \\
0.4 & 2.4 & 0.563 \\
0.3 & 1.4 & 0.314 \\
0.2 & 0.6 & 0.095 \\
\hline
\end{tabular}}
\hspace{0.1\textwidth}
\subfigure[$n=400$, $\varepsilon=0.2$, $m\approx 8000$]{
\label{tab:small400}
\begin{tabular}{|c|c|c|}
\hline
$p$ & $\frac{2mp^2}{n}$ & $\rho_1$ \\
\hline
0.8 & 26 & 0.960 \\
0.4 & 6.4 & 0.817 \\
0.3 & 3.6 & 0.643 \\
0.2 & 1.6 & 0.282 \\
0.1 & 0.4 & 0.145 \\
\hline
\end{tabular}}
\end{center}
\caption{Correlations between the top eigenvector of $H$ and the vector of true angles for different values of $p$ in the small-world $S^2$ model.} \label{tab:small}
\end{table}

The comparison between the eigenvector and SDP methods (as well as the least squares method of Section \ref{sec:intro}) is summarized in Table \ref{tab:sdp-eig2} showing the numerical correlations for $n=200$, $\varepsilon=0.3$ (number of edges $m\approx 3000$) and for different values of $p$. Although the SDP is slightly more accurate, the eigenvector method runs faster.

\begin{table}[h]
\begin{center}
\begin{tabular}{|c||c|c|c|c|}
\hline
$p$ & $\rho_{lsqr}$ & $\rho_{eig}$ & $\rho_{sdp}$ & $\operatorname{rank} \Theta$ \\
\hline
1   & 1 & 1 & 1 & 1\\
0.7 & 0.787 & 0.977 & 0.986 & 1 \\
0.4 & 0.046 & 0.839 & 0.893 & 3 \\
0.3 & 0.103 & 0.560 & 0.767 & 3 \\
0.2 & 0.227 & 0.314 & 0.308 & 4 \\
0.15& 0.091 & 0.114 & 0.102 & 5 \\
\hline
\end{tabular}
\end{center}
\caption{Comparison between the correlations obtained by the eigenvector method $\rho_{eig}$, by the SDP method $\rho_{sdp}$ and by the least squares method $\rho_{lsqr}$ for different values of $p$ (small world graph on $S^2$, $n=200$, $\varepsilon=0.3$, $m\approx 3000$). The SDP tends to find low-rank matrices despite the fact that the rank-one constraint on $\Theta$ is not included in the SDP.  The rightmost column gives the rank of the $\Theta$ matrices that were found by the SDP. To solve the SDP (\ref{sdp1a})-(\ref{sdp1c}) we used \texttt{SDPLR}, a package for solving large-scale SDP problems \cite{SDPLR}. The least squares solution was obtained using MATLAB's \textsf{lsqr} function. As expected, the least squares method yields poor correlations compared to the eigenvector and the SDP methods.} \label{tab:sdp-eig2}
\end{table}

\section{Information Theoretic Analysis}
\label{sec:information-theory}
The optimal solution to the angular synchronization problem can be considered as the set of angles that maximizes the log-likelihood. Unfortunately, the log-likelihood is a non-convex function and the maximum likelihood cannot be found in a polynomial time. Both the eigenvector method and the SDP method are polynomial-time relaxations of the maximum log-likelihood problem. In the previous section we showed that the eigenvector method fails to recover the true angles when $p$ is below the threshold probability $p_c^{\textit{eig}}$. It is clear that even the maximum likelihood solution would fail to recover the correct set of angles below some (perhaps lower) threshold. It is therefore natural to ask if the threshold value of the polynomial eigenvector method gets close to the optimal threshold value of the exponential-time maximum likelihood exhaustive search. In this section we provide a positive answer to this question using the information theoretic Shannon bound \cite{Cover}. Specifically, we show that the threshold probability for the eigenvector method is asymptotically  larger by just a multiplicative factor compared to the threshold probability of the optimal recovery algorithm. The multiplicative factor is a function of the angular discretization resolution, but not a function of $n$ and $m$. The eigenvector method becomes less optimal as the discretization resolution improves.

We start the analysis by recalling that from the information theoretic point of view, the uncertainty in the values of the angles is measured by their entropy.
The noisy offset measurements carry some bits of information on the angle values, therefore decreasing their uncertainty, which is measured by the conditional entropy that we need to estimate.

The angles $\theta_1,\ldots,\theta_n$ can take any real value in the interval $[0,2\pi)$. However, an infinite number of bits is required to describe real numbers, and so we cannot hope to determine the angles with an arbitrary precision. Moreover, the offset measurements are often also discretized. We therefore seek to determine the angles only up to some discretization precision $\frac{2\pi}{L}$, where $L$ is the number of subintervals of $[0,2\pi)$ obtained by dividing the unit circle is into $L$ equally sized pieces.

Before observing any of the offset measurements, the angles are uniformly distributed on $\{0,1,\ldots,L-1\}$, that is, each of them falls with equal probability $1/L$ to any of the $L$ subintervals. It follows that the entropy of the $i$'th angle $\theta_i$ is given by
\begin{equation}
\label{theta1}
H(\theta_i) = -\sum_{l=0}^{L-1} \frac{1}{L}\log_2 \frac{1}{L} = \log_2 L,\quad \mbox{for } i=1,2,\ldots,n.
\end{equation}
We denote by $\mb{\theta}^n = (\theta_1,\ldots,\theta_n)$ the vector of angles.
Since $\theta_1,\ldots,\theta_n$ are independent (the orientations of the molecules are random), their joint entropy $H(\mb{\theta}^n)$ is given by
\begin{equation}
\label{entropy-angles}
H(\mb{\theta}^n) = \sum_{i=1}^n H(\theta_i) = n\log_2 L,
\end{equation}
reflecting the fact that the configuration space is of size $L^n = 2^{n\log_2 L}$.

Let $\delta_{ij}$ be the random variable for the outcome of the noisy offset measurement of $\theta_i$ and $\theta_j$. The random variable $\delta_{ij}$ is also discretized and takes values in $\{0,1,\ldots,L-1\}$. We denote by $\mb{\delta}^m = (\delta_{i_1j_1}$,\ldots,$\delta_{i_mj_m})$ the vector of all offset measurements. Conditioned on the values of $\theta_i$ and $\theta_j$, the random variable $\delta_{ij}$ has the following conditional probability distribution
\begin{equation}
\label{eq:offset-dist}
\Pr\{\delta_{ij} \,|\,\theta_i,\theta_j \} = \left\{\begin{array}{cc} \frac{1-p}{L} & \delta_{ij} \neq \theta_i-\theta_j \mod L, \\ p + \frac{1-p}{L} & \delta_{ij}=\theta_i-\theta_j \mod L,
\end{array} \right.
\end{equation}
because with probability $1-p$ the measurement $\delta_{ij}$ is an outlier that takes each of the $L$ possibilities with equal probability $\frac{1}{L}$, and with probability $p$ it is a good measurement that equals $\theta_i-\theta_j$.
It follows that the conditional entropy $H(\delta_{ij} | \theta_i,\theta_j)$ is
\begin{equation}
\label{entropy-offset}
H(\delta_{ij} | \theta_i,\theta_j) = -(L-1)\frac{1-p}{L}\log_2\frac{1-p}{L} - \left(p + \frac{1-p}{L}\right)\log_2 \left(p + \frac{1-p}{L}\right).
\end{equation}
We denote this entropy by $H(L,p)$ and its deviation from $\log_2 L$ by $I(L,p)$, that is,
\begin{equation}
\label{H11}
H(L,p) \equiv -(L-1)\frac{1-p}{L}\log_2\frac{1-p}{L} - \left(p + \frac{1-p}{L}\right)\log_2 \left(p + \frac{1-p}{L}\right).
\end{equation}
and
\begin{equation}
\label{I11}
I(L,p) \equiv  \log_2 L - H(L,p).
\end{equation}
Without conditioning, the random variable $\delta_{ij}$ is uniformly distributed on $\{0,\ldots,L-1\}$ and has entropy
\begin{equation}
H(\delta_{ij}) = \log_2 L.
\end{equation}
It follows that the mutual information $I(\delta_{ij}; \theta_i,\theta_j)$ between the offset measurement $\delta_{ij}$ and the angle values $\theta_i$ and $\theta_j$ is
\begin{equation}
I(\delta_{ij}; \theta_i,\theta_j) = H(\delta_{ij})-H(\delta_{ij}|\theta_i,\theta_j) = \log_2 L - H(L,p) = I(L,p).
\end{equation}
This mutual information measures the reduction in the uncertainty of the random variable $\delta_{ij}$ from knowledge of $\theta_i$ and $\theta_j$. Due to the symmetry of the mutual information,
\begin{equation}
I(\delta_{ij}; \theta_i,\theta_j) = H(\delta_{ij})-H(\delta_{ij}|\theta_i,\theta_j) = H(\theta_i,\theta_j)-H(\theta_i,\theta_j|\delta_{ij}),
\end{equation}
the mutual information is also the reduction in uncertainty of the angles $\theta_i$ and $\theta_j$ given the noisy measurement of their offset $\delta_{ij}$. Thus,
\begin{equation}
H(\theta_i,\theta_j | \delta_{ij}) = H(\theta_i,\theta_j) - I(\delta_{ij}; \theta_i,\theta_j).
\end{equation}
Similarly, given all $m$ offset measurements $\mb{\delta}^m$, the uncertainty in $\mb{\theta}^n$ is given by
\begin{equation}
\label{cond-entropy}
H(\mb{\theta}^n | \mb{\delta}^m) = H(\mb{\theta}^n) - I(\mb{\delta}^m ; \mb{\theta}^n),
\end{equation}
with
\begin{equation}
\label{mutual-information}
I(\mb{\delta}^m ; \mb{\theta}^n) = H(\mb{\delta}^m) - H(\mb{\delta}^m|\mb{\theta}^n).
\end{equation}
A simple upper bound for this mutual information is obtained by explicit evaluation of the conditional entropy $H(\mb{\delta}^m|\mb{\theta}^n)$ combined with a simple upper bound on the joint entropy term $H(\mb{\delta}^m)$. First, note that given the values of $\theta_1,\ldots,\theta_n$, the offsets become independent random variables. That is, knowledge of $\delta_{i_1j_1}$ (given $\theta_{i_1},\theta_{j_1}$) does not give any new information on the value of $\delta_{i_2j_2}$ (given $\theta_{i_2},\theta_{j_2}$). The conditional probability distribution of the offsets is completely determined by (\ref{eq:offset-dist}), and the conditional entropy is therefore the sum of $m$ identical entropies of the form (\ref{entropy-offset})
\begin{equation}
\label{entropy-2}
H(\mb{\delta}^m|\mb{\theta}^n) = m H(L,p).
\end{equation}
Next, bounding the joint entropy $H(\mb{\delta}^m)$ by the logarithm of its configuration space size $L^m$ yields
\begin{equation}
\label{entropy-bound}
H(\mb{\delta}^m) \leq m\log_2 L.
\end{equation}
Note that this simple upper bound ignores the dependencies among the offsets which we know to exist, as implied, for example, by the triplet consistency relation (\ref{consistent}). As such, (\ref{entropy-bound}) is certainly not a tight bound, but still good enough to prove our claim about the nearly optimal performance of the eigenvector method.

Plugging (\ref{entropy-2}) and (\ref{entropy-bound}) in (\ref{mutual-information}) yields the desired upper bound on the mutual information
\begin{equation}
\label{mutual-bound}
I(\mb{\delta}^m ; \mb{\theta}^n) \leq m\log_2 L - m H(L,p) = m I(L,p).
\end{equation}
Now, substituting the bound (\ref{mutual-bound}) and the equality (\ref{entropy-angles}) in (\ref{cond-entropy}) gives a lower bound for the conditional entropy
\begin{equation}
\label{bbb}
H(\mb{\theta}^n | \mb{\delta}^m) \geq n\log_2 L - mI(L,p).
\end{equation}
We may interpret this bound in the following way. Before seeing any offset measurement the entropy of the angles is $n\log_2 L$, and each of the $m$ offset measurements can decrease the conditional entropy by at most $I(L,p)$, the information that it carries.

The bound (\ref{bbb}) demonstrates, for example, that for fixed $n$, $p$ and $L$, the conditional entropy is bounded from below by a linear decreasing function of $m$. It follows that unless $m$ is large enough, the uncertainty in the angles would be too large. Information theory says that a successful recovery of all $\theta_1,\ldots,\theta_n$ is possible only when their uncertainty, as expressed by the conditional entropy, is small enough.
The last statement can be made precise by Fano's inequality and Wolfowitz' converse, also known as the weak and strong converse theorems to the coding theorem that provide a lower bound for the probability of the error probability in terms of the conditional entropy, see, e.g., \cite[Chapter 8.9, pages 204-207]{Cover} and \cite[Chapter 5.8, pages 173-176]{Gallager}.

In the language of coding, we may think of $\mb{\theta}^n$ as a codeword that we are trying to decode from the noisy vector of offsets $\mb{\delta}^m$ which is probabilistically related to $\mb{\theta}^n$. The codeword $\mb{\theta}^n$ is originally uniformly distributed on $\{1,2,\ldots,2^{n\log_2 L}\}$ and from $\mb{\delta}^m$ we estimate $\mb{\theta}^n$ as one of the $2^{n\log_2 L}$ possibilities. Let the estimate be ${\mb{\hat{\theta}}}^n$ and define the probability of error as $P_e = \Pr\{\mb{\hat{\theta}}^n \neq \mb{\theta}^n\}$. Fano's inequality \cite[Lemma 8.9.1, page 205]{Cover} gives the following lower bound on the error probability
\begin{equation}
\label{combine}
H(\mb{\theta}^n | \mb{\delta}^m) \leq 1 + P_e n\log_2 L.
\end{equation}
Combining (\ref{combine}) with the lower bound for the conditional entropy (\ref{bbb}) we obtain a weak lower bound on the error probability
\begin{equation}
\label{fano}
P_e \geq 1 - \frac{m}{n}\frac{I(L,p)}{\log_2 L} - \frac{1}{n\log_2 L}.
\end{equation}
This lower bound for the probability of error is applicable to all decoding algorithms, not just for the eigenvector method. For large $n$, we see that for any $\beta < 1$,
\begin{equation}
\label{fano2}
\frac{m}{n}\frac{I(L,p)}{\log_2 L} < \beta  \Longrightarrow P_e \geq 1 - \beta + o(1).
\end{equation}
We are mainly interested in the limit $m,n \to \infty$ and $p\to 0$ with $L$ being fixed. The Taylor expansion of $I(L,p)$ (given by (\ref{H11})-(\ref{I11})) near $p=0$ reads
\begin{equation}
\label{taylor}
I(L,p) = \frac{1}{2} (L-1)p^2 + O(p^3).
\end{equation}
Combining (\ref{fano2}) and (\ref{taylor}) we obtain that
\begin{equation}
\label{fano3}
p = \sqrt {\frac{n}{m}\frac{2\log_2 L}{(L-1)}\beta} \Longrightarrow P_e \geq 1-\beta + o(1), \quad \mbox{as } n,m\to \infty,\; n/m\to 0.
\end{equation}
Note that $n/m\to 0$, because $m\geq n\log n$ in order to ensure with high probability the connectivity of the measurement graph $G$.
The bound (\ref{fano3}) was derived using the weak converse theorem (Fano's inequality). It is also possible to show that the probability of error goes exponentially to 1 (using the Wolfowitz' converse and Chernoff bound, see \cite[Theorem 5.8.5, pages 173-176]{Gallager}).

The above discussion shows that there does not exist a decoding algorithm with a small probability for the error for values of $p$ below the threshold probability $p_c^{\textit{inf}}$ given by
\begin{equation}
\label{pc-opt}
p_c^{\textit{inf}} = \sqrt{\frac{n}{m}\frac{2\log_2 L}{L-1}}.
\end{equation}
Note that for $L=2$, the threshold probability $p_c^{\textit{eig}} = \frac{1}{\sqrt{n}}$ of the eigenvector method in the complete graph case for which $m={n \choose 2}$ is 2 times smaller than $p_c^{inf}$. This is not a violation of information theory, because the fact that the top eigenvector has a non-trivial correlation with the vector of true angles does not mean that all angles are recovered correctly by the eigenvector. The fact that the eigenvector method gives non-trivial correlations below the information theoretic bound is just another evidence to its effectiveness.

We turn to shed some light on why it is possible to partially recover the angles below the information theoretic bound.
The main issue here is that it is perhaps too harsh to measure the success of the decoding algorithm by $P_e=\Pr\{\mb{\hat{\theta}}^n \neq \mb{{\theta}}^n\}$. For example, when the decoding algorithm decodes $999$ angles out of $n=1000$ correctly while making just a single mistake, we still count it as a failure. It may be more natural to consider the probability of error in the estimation of the individual angles. We proceed to show that this measure of error leads to a threshold probability which is smaller than (\ref{pc-opt}) by just a constant factor.

Let $P_e^{(1)}=\Pr\{\hat{\theta}_1 \neq \theta_1 \}$ be the probability of error in the estimation of $\theta_1$. Again, we want to use Fano's inequality to bound the probability of the error by bounding the conditional entropy. A simple lower bound to the conditional entropy $H(\theta_1 | \mb{\delta}^m)$ is obtained by conditioning on the remaining $n-1$ angles
\begin{equation}
\label{theta-bound}
H(\theta_1 | \mb{\delta}^m) \geq H(\theta_1 | \mb{\delta}^m,\theta_2,\theta_3,\ldots,\theta_n).
\end{equation}
Suppose that there are $d_1$ noisy offset measurements of the form $\theta_1-\theta_j$, that is, $d_1$ is the degree of node 1 in the measurement graph $G$. Let the neighbors of node 1 be $j_1,j_2,\ldots,j_{d_1}$ with corresponding offset measurements $\delta_{1j_1},\ldots,\delta_{1j_{d_1}}$. Given the values of all other angles $\theta_2,\ldots,\theta_n$, and in particular the values of $\theta_{j_1},\ldots,\theta_{j_{d_1}}$, these $d_1$ equations become noisy equations for the single variable $\theta_1$. We denote these transformed equations for $\theta_1$ alone by $\tilde{\delta}_{1},\ldots,\tilde{\delta}_{d_1}$. All other $m-d_1$ equations do not involve $\theta_1$ and therefore do not carry any information on its value. It follows that
\begin{equation}
\label{tilde0}
H(\theta_1 | \mb{\delta}^m,\theta_2,\theta_3,\ldots,\theta_n) = H(\theta_1 | \tilde{\delta}_{1},\ldots,\tilde{\delta}_{d_1}).
\end{equation}
We have
\begin{equation}
\label{tilde1}
H(\tilde{\delta}_1,\ldots,\tilde{\delta}_{d_1} | \theta_1 ) = d_1 H(L,p),
\end{equation}
because given $\theta_1$ these $d_1$ equations are i.i.d random variables with entropy $H(L,p)$. Also, a simple upper bound on the $d_1$ equations (without conditioning) is given by
\begin{equation}
\label{tilde2}
H(\tilde{\delta}_1,\ldots,\tilde{\delta}_{d_1} ) \leq d_1 \log_2 L,
\end{equation}
ignoring possible dependencies among the outcomes. From (\ref{tilde1})-(\ref{tilde2}) we get an upper bound for the mutual information between $\theta_1$ and the transformed equations
\begin{equation}
\label{tilde3}
I(\theta_1 ; \tilde{\delta}_1,\ldots,\tilde{\delta}_{d_1}) \leq d_1 \left[\log_2 L - H(L,p)\right] = d_1 I(L,p).
\end{equation}
Combining (\ref{theta-bound}),(\ref{tilde0}) (\ref{tilde3}) and (\ref{theta1}) we get
\begin{eqnarray}
H(\theta_1 | \mb{\delta}^m) &\geq& H(\theta_1 | \mb{\delta}^m,\theta_2,\theta_3,\ldots,\theta_n) \label{t-bound}\\
&=& H(\theta_1 | \tilde{\delta}_{1},\ldots,\tilde{\delta}_{d_1}) \nonumber \\
&=& H(\theta_1) - I(\theta_1 ; \tilde{\delta}_{1},\ldots,\tilde{\delta}_{d_1}) \nonumber \\
&\geq& \log_2 L - d_1 I(L,p).\nonumber
\end{eqnarray}
This lower bound on the conditional entropy translates, via Fano's inequality, to a lower bound on the probability of error $P_e^{(1)}$, and it follows that
\begin{equation}
d_1 I(L,p) > \log_2 L
\end{equation}
is a necessary condition for having a small $P_e^{(1)}$. Similarly, the condition for a small probability of error in decoding $\theta_i$ is
\begin{equation}
\label{trim1}
d_i I(L,p) > \log_2 L,
\end{equation}
where $d_i$ is the degree of vertex $i$ in the measurement graph. This condition suggests that we should have more success in decoding angles of high degree. The average degree $\bar{d}$ in a graph with $n$ vertices and $m$ edges is $\bar{d} = \frac{2m}{n}$. The condition for successful decoding of angles with degree $\bar{d}$ is
\begin{equation}
\label{trim2}
\frac{2m}{n} I(L,p) > \log_2 L.
\end{equation}
In particular, this would be the condition for all vertices in a regular graph, or in a graph whose degree distribution is concentrated near $\bar{d}$.

Substituting the Taylor expansion (\ref{taylor}) into (\ref{trim2}) results in the condition
\begin{equation}
p > \sqrt{\frac{n}{m}\frac{\log_2 L}{L-1}}.
\end{equation}
This means that successful decoding of the individual angles may be possible already for $p > p_c^{\textit{ind}}$, where
\begin{equation}
\label{p-ind}
p_c^{\textit{ind}} = \sqrt{\frac{n}{m}\frac{\log_2 L}{L-1}},
\end{equation}
but the estimation of the individual angles must contain some error when $p < p_c^{\textit{ind}}$. Note that $p_c^{\textit{ind}} < p_c^{\textit{inf}}$, so while for $p$ values between $p_c^{\textit{ind}}$ and $p_c^{\textit{inf}}$ it is impossible to successfully decode all angles, it may still be possible to decode some angles.

In the complete graph case, comparing the threshold probability of the eigenvector method $p_c^{\textit{eig}} = \frac{1}{\sqrt{n}}$ given by (\ref{pn}) and the information theoretic threshold probability $p_c^{\textit{ind}}$ (\ref{p-ind}) below which no algorithm can successfully recover individual angles, we find that their ratio is asymptotically independent of $n$ and $m$:
\begin{equation}
\label{p-eig-ind}
\frac{p_c^{\textit{eig}}}{p_c^{\textit{ind}}} = \sqrt{\frac{L-1}{2\log_2 L}} + o(1).
\end{equation}
Note that the threshold probability $p_c^{\textit{eig}}$ is smaller than $p_c^{\textit{ind}}$ for $L\leq 6$. Thus, we may regard the eigenvector method as a very successful recovery algorithm for offset equations with a small modulo $L$.

For $L \geq 7$, equation (\ref{p-eig-ind}) implies a gap between the threshold probabilities $p_c^{\textit{eig}}$ and $p_c^{\textit{ind}}$, suggesting that the exhaustive exponential search for the maximum likelihood would perform better than the polynomial time eigenvector method. Note, however, that the gap would be significant only for very large values of $L$ that correspond to very fine angular resolutions. For example, even for $L=100$ the threshold probability of the eigenvector method would only be $\sqrt{\frac{99}{2\log_2 100}} \approx 2.73$ times larger than that of the maximum likelihood. The exponential complexity of $O(mL^n)$ of the exhaustive search for the maximum likelihood makes it impractical even for moderate-scale problems. On the other hand, the eigenvector method has a polynomial running time and it can handle large scale problems with relative ease.

\section{Connection with \textsc{Max-2-Lin mod} $L$ and Unique Games}
\label{sec:max2lin}

The angular synchronization problem is related to the combinatorial optimization problem \textsc{Max-2-Lin mod} $L$ for maximizing the number of satisfied linear equations mod $L$ with exactly 2 variables in each equation, because the discretized offset equations $\theta_i-\theta_j = \delta_{ij} \mod L$ are exactly of this form. \textsc{Max-2-Lin mod} $L$ is a problem mainly studied in theoretical computer science, where we prefer using the notation ``\textsc{mod $L$}" instead of the more common ``\textsc{mod $p$}", to avoid confusion between the size of the modulus and the proportion of good measurements.

Note that a random assignment of the angles would satisfy a $\frac{1}{L}$ fraction of the offset equations.
Andersson, Engebretsen, and H{\aa}stad \cite{Hastad} considered SDP based algorithms
for \textsc{Max-2-Lin mod} $L$, and showed that they could obtain an $\frac{1}{L}\left(1+\kappa(L)\right)$-approximation algorithm, where $\kappa(L) > 0$ is a constant that depends on $L$. In particular, they gave a very weak proven performance guarantee of $\frac{1}{L}(1 + 10^{-8})$, though they concluded that it is most likely that their bounds can be improved significantly. Moreover, for $L=3$ they numerically find the approximation ratio to be $\frac{1}{1.27} \approx 0.79$, and later Goemans and Williamson \cite{GoemansMAX3CUT} proved a 0.793733-approximation. The SDP based algorithms in \cite{Hastad} are similar in their formulation to the SDP based algorithm of Frieze and Jerrum for \textsc{Max}-k-\textsc{Cut} \cite{Max-k-Cut}, but with a different rounding procedure. In these SDP models, $L$ vectors are assigned to each of the $n$ angle variables, so that the total number of vectors is $nL$. The resulting $nL \times nL$ matrix of inner products is required to be semidefinite positive, along with another set of $O(n^2L^2)$ linear and inequality constraints. Due to the large size of the inner product matrix and the large number of constraints, our numerical experiments with these SDP models were limited to relatively small size problems (such as $n=20$ and $L=7$) from which it was difficult to get a good understanding of their performance. In the small scale problems that we did manage to test, we did not find any supporting evidence that these SDP algorithms perform consistently better than the eigenvector method, despite their extensive running times and memory requirements. For our SDP experiments we used the software \texttt{SDPT3} \cite{SDPT31,SDPT32} and \texttt{SDPLR} \cite{SDPLR} in MATLAB. In \cite{Hastad} it is also shown that it is NP-hard to approximate \textsc{Max-2-Lin mod} $L$ within a constant ratio, independent of $L$. Thus, we should expect an $L$-dependent gap similar to (\ref{p-eig-ind}) for any polynomial time algorithm, not just for the eigenvector method.

\textsc{Max-2-Lin} is an instance of what is known as unique games \cite{Feige}, described below. One distinguishing feature of the offset equations is that every constraint
corresponds to a bijection between the values of the associated variables. That is, for every possible value
of $\theta_i$, there is a unique value of $\theta_j$ that satisfies the constraint $\theta_i-\theta_j=\delta_{ij}$. Unique
games are systems of constraints, a generalization of the offset equations, that have
this uniqueness property, so that every constraint corresponds to some permutation.

As in the setting of offset equations, instances of unique games where all constraints are satisfiable
are easy to handle. Given an instance where $1-\varepsilon$ fraction of constraints are satisfiable, the
Unique Games Conjecture (UGC) of Khot \cite{Khot} says that it is hard to satisfy even a $\gamma>0$ fraction of the
constraints. The UGC has been shown to imply a number of inapproximability
results for fundamental problems that seem difficult to obtain by more standard complexity
assumptions. Note that in our angular synchronization problem the fraction of constraints that are satisfiable is $1-\varepsilon = p + \frac{1-p}{L}$.

Charikar, Makarychev and Makarychev \cite{Charikar} presented improved approximation algorithms for unique games. For instances with domain size
$L$ where the optimal solution satisfies $1-\varepsilon$ fraction of all constraints, their algorithms satisfy
roughly $L^{-\varepsilon/(2-\varepsilon)}$ and $1 - O(\sqrt{\varepsilon \log L})$ fraction of all constraints. Their algorithms are based on SDP, also with an underlying  inner products matrix of size $nL \times nL$, but their constraints and rounding procedure are different than those of \cite{Hastad}. Given the results of \cite{Odonell}, the algorithms in \cite{Charikar} are near optimal if the UGC is true, that is, any improvement (beyond low order terms) would refute the
conjecture. We have not tested their SDP based algorithm in practice, because, like the SDP of \cite{Hastad} it is also expected to be limited to relatively small scale problems.

\section{Summary and Further Applications}
\label{sec:sync}

In this paper we presented an eigenvector method and an SDP approach for solving the angular synchronization problem. We used random matrix theory to prove that the eigenvector method finds an accurate estimate for the angles even in the presence of a large number of outlier measurements.

The idea of synchronization by eigenvectors can be applied to other problems exhibiting a group structure and noisy measurements of ratios of group elements. In this paper we specialized the synchronization problem over the group SO(2). In the general case we may consider a group $G$ other than SO(2) for which we have good and bad measurements $g_{ij}$ of ratios between group elements
\begin{equation}
g_{ij} = g_i {g_j}^{-1}, \quad g_i,{g_j} \in G.
\end{equation}
For example, in the general case, the triplet consistency relation (\ref{consistent}) simply reads
\begin{equation}
g_{ij}g_{jk}g_{ki} = g_i {g_j}^{-1} g_j {g_k}^{-1} g_k {g_i}^{-1} = e,
\end{equation}
where $e$ is the identity element of $G$.

Whenever the group $G$ is compact and has a complex or real representation
(for example, the rotation group SO(3) has a real representation using $3\times 3$ rotation matrices),
we may construct an Hermitian matrix that is a matrix of matrices: the $ij$ element is either the matrix representation of the measurement $g_{ij}$ or the zero matrix if there is no direct measurement for the ratio of $g_i$ and $g_j$. Once the matrix is formed, once can look for its top eigenvectors (or SDP) and estimate the group elements from them.

In some cases the eigenvector and the SDP methods can be applied even when there is only partial information for the group ratios. This problem arises naturally in the determination of the three-dimensional structure of a macromolecule in cryo-electron microscopy \cite{Fred1}. In \cite{common-lines} we show that the common lines between projection images give partial information for the group ratios between elements in SO(3) that can be estimated accurately using the eigenvector and SDP methods.
In \cite{class-averaging} we explore the close connection between the angular synchronization problem and the class averaging problem in cryo-electron microscopy \cite{Fred1}. Other possible applications of the synchronization problem over SO(3) include the distance geometry problem in NMR spectroscopy \cite{DC1,DC4} and the localization of sensor networks \cite{Ye2006,AmitPNAS}.

The eigenvector method can also be applied to non-compact groups that can be ``compactified". For example, consider the group of real numbers $\mathbb{R}$ with addition.
One may consider the synchronization problem of clocks that measure noisy time differences of the form
\begin{equation}
\label{clocks}
t_i - t_j = t_{ij}, \quad t_i,t_j \in \mathbb{R}.
\end{equation}
We compactify the group $\mathbb{R}$ by mapping it to the unit circle $t \mapsto e^{\imath \omega t}$, where $\omega \in \mathbb{R}$ is a parameter
to be chosen not too small and not too large, as we now explain.
There may be two kinds of measurement errors in (\ref{clocks}). The first kind of error is a small discretization error (e.g., a small Gaussian noise)
of typical size $\Delta$. The second type of error is a large error that can be regarded as an outlier. For example, in some practical application an error of size
$10\Delta$ may be considered as an outlier. We therefore want $\omega$ to satisfy $\omega \gg (1/10)\Delta^{-1}$ (not too small) and
$\omega \ll  \Delta^{-1}$ (not too large), so that when constructing the matrix
\begin{equation}
\label{H-clocks}
H_{ij} = \left\{\begin{array}{cc}
 e^{\imath \omega t_{ij}} & \{i,j\}\in E, \\
 0 & \{i,j\} \not\in E,
\end{array}\right.
\end{equation}
each good equation will contribute approximately 1, while the contribution of the bad equations will be uniformly distributed on the unit circle.
One may even try several different values for the ``frequency" $\omega$ in analogy to the Fourier transform.
An overdetermined linear system of the form (\ref{clocks}) can also be solved using least squares, which is also the maximum likelihood estimator if the measurement errors are Gaussian.
However, in the many outliers model, the contribution of outlier equations will dominate the sum of squares error. For example, each outlier equation with error $10\Delta$
contributes to the sum of squares error the same as 100 good equations with error $\Delta$. The compactification of the group combined with the eigenvector method
has the appealing effect of reducing the impact of the outlier equations.
This may open the way for the eigenvector method based on (\ref{H-clocks}) to be useful for the surface reconstruction problems in computer vision \cite{FrankotChellappa,AgrawalWhatRange} and optics
\cite{Koby2001} in which current methods succeed only in the presence of a limited number of outliers.

\section*{Acknowledgments}
The author would like to thank Yoel Shkolniskly, Fred Sigworth and Ronald Coifman for many stimulating discussions regarding the cryo-electron microscopy problem;
Boaz Barak for references to the vast literature on \textsc{Max-2-Lin mod} $L$ and unique games;
Amir Bennatan for pointers to the weak and strong converse theorems to the coding theorem;
Robert Ghrist and Michael Robinson for valuable discussions at UPenn and for the reference to \cite{GiridharKumar};
and Steven (Shlomo) Gortler, Yosi Keller and Ben Sonday for reviewing an earlier version of the manuscript and for their helpful suggestions.

The project described was supported by Award Number DMS-0914892 from the NSF, by Award Number FA9550-09-1-0551 from AFOSR, and by Award Number R01GM090200 from the National Institute Of
General Medical Sciences. The content is solely the responsibility of the authors and does not
necessarily represent the official views of the National Institute Of General Medical Sciences or the
National Institutes of Health.


\begin{thebibliography}{99}

\bibitem{AgrawalWhatRange}
Agrawal, A.~K., Raskar, R., and Chellappa, R. (2006) What is the range of surface
  reconstructions from a gradient field? in \emph{Computer Vision -- ECCV
  2006: 9th European Conference on Computer Vision, Graz, Austria, May 7-13,
  2006, Proceedings, Part IV (Lecture Notes in Computer Science)}, pp.
  578--591.

\bibitem{AlonVu}
Alon, N., Krivelevich, M., and Vu, V. H. (2002) On the concentration of eigenvalues of random symmetric matrices,
{\em Israel Journal of Mathematics} {\bf 131} (1) pp. 259--267.

\bibitem{Hastad}
Andersson, G., Engebretsen, L., and H{\aa}stad, J. (1999)
A new way to use semidefinite programming with applications to linear equations mod $p$.
{\em Proceedings 10th annual ACM-SIAM symposium on Discrete algorithms}, pp. 41--50.

\bibitem{Ye2006} {Biswas, P.}, {Liang, T. C.}, {Toh, K.C.}, {Wang, T. C.}, and {Ye, Y.} (2006)
\newblock Semidefinite programming approaches for sensor network localization with noisy distance measurements.
\newblock {\em IEEE Transactions on Automation Science and Engineering},
  {\bf 3}(4):360--371.

\bibitem{SDPLR}
Burer, S. and Monteiro, R.D.C. (2003)
A Nonlinear Programming Algorithm for Solving Semidefinite Programs Via Low-Rank Factorization.
{\em Mathematical Programming (series B)}, {\bf 95} (2):329--357.
%

\bibitem{Charikar}
Charikar, M., Makarychev, K., and Makarychev, Y. (2006)
Near-Optimal Algorithms for Unique Games.
{\em Proceedings 38th annual ACM symposium on Theory of computing}, pp. 205--214.


\bibitem{Coifman2006a}
{Coifman, R. R.} and {Lafon, S.} (2006) Diffusion maps. \emph{Applied and
  Computational Harmonic Analysis}, {\bf 21} (1), pp. 5--30.

\bibitem{Cover}
Cover, T. M. and Thomas, J. A. (1991) {\em Elements of Information Theory},
Wiley, New York.


\bibitem{ErdosRenyi}
Erd\H{o}s, P. and R\'enyi, A. (1959)
On random graphs. {\em Publicationes Mathematicae} {\bf 6}, pp. 290--297.

\bibitem{Feige}
Feige, U., and Lov{\'a}sz, L. (1992)
Two-prover one round proof systems: Their power and their problems.
{\em In Proceedings of the 24th ACM Symposium on Theory of Computing}, pp. 733--741.

\bibitem{PecheFeral}
F\'eral, D.  and P\'ech\'e, S. (2007)
The Largest Eigenvalue of Rank One Deformation of Large Wigner Matrices, {\em Communications in Mathematical Physics}
{\bf 272} (1): 185--228.

\bibitem{Fred1} {Frank, J.} (2006) {\em Three-Dimensional Electron Microscopy of Macromolecular Assemblies: Visualization of Biological Molecules in Their Native State},
Oxford.

\bibitem{FrankotChellappa}
Frankot, R.~T. and Chellappa, R. (1988) A method for enforcing integrability in shape
  from shading algorithms, {\em IEEE Transactions on Pattern Analysis and Machine Intelligence}, {\bf 10} (4): 439-451.

\bibitem{Max-k-Cut}
Frieze A., and Jerrum, M. (1997)
Improved Approximation Algorithms for MAX k-CUT and MAX BISECTION.
{\em Algorithmica} {\bf 18} (1), pp. 67--81.

\bibitem{Furedi}
F\"uredi, Z., and Koml\'{o}s, J. (1981)
The eigenvalues of random symmetric matrices.
{\em Combinatorica}, {\bf 1}, pp. 233--241.

\bibitem{Gallager}
Gallager, R. G. (1968) {\em Information Theory and Reliable Communication}, Wiley, New York.

\bibitem{GiridharKumar}
Giridhar, A. and Kumar, P.R. (2006)
Distributed Clock Synchronization over Wireless Networks: Algorithms and Analysis,
{\em 45th IEEE Conference on Decision and Control 2006}, pp. 4915--4920.

\bibitem{GoemansWilliamson}
Goemans, M. X. and Williamson, D. P. (1995)
Improved approximation algorithms for maximum cut and satisfiability problems using semidefinite programming.
{\em Journal of the ACM (JACM)} {\bf 42} (6), pp.~1115--1145.

\bibitem{GoemansMAX3CUT}
Goemans, M. X. and Williamson, D. P. (2001)
Approximation algorithms for \textsc{Max-3-Cut} and other problems via complex semidefinite programming.
{\em Proceedings 33rd annual ACM symposium on Theory of computing}, pp. 443--452.
%
%

\bibitem{Griffiths}
Griffiths, D. J. (1994)
{\em Introduction to Quantum Mechanics}, Prentice Hall, NJ, 416 pages.

\bibitem{DC4} Havel, T.~F., and Wuthrich, K. (1985). An evaluation of the combined use of nuclear magnetic resonance and distance geometry for the determination of
    protein conformation in solution. {\em J Mol Biol} {\bf 182}, 281-294.


\bibitem{HornMatrix}
Horn, R. A., and Johnson, C. R. (1990)
{\em Matrix Analysis}. Cambridge University Press, 575 pages.

\bibitem{Karp}
Karp, R., Elson, J., Estrin, D. and Shenker, S. (2003) Optimal and
global time synchronization in sensornets, {\em Technical Report}, Center
for Embedded Networked Sensing, University of California,
Los Angeles.

\bibitem{Khorunzhiy2003}
Khorunzhiy, O. (2003) Rooted trees and moments of large sparse
random matrices, {\em Discrete Mathematics and Theoretical Computer Science AC}, pp. 145--154.

\bibitem{Khorunzhy2001}
Khorunzhy, A. (2001)
Sparse Random Matrices: Spectral Edge and Statistics of Rooted Trees,
{\em Advances in Applied Probability} {\bf 33} (1) pp. 124--140.


\bibitem{Khot}
Khot, S. (2002)
On the power of unique 2-prover 1-round games.
In {\em Proceedings of the ACM
Symposium on the Theory of Computing}, {\bf 34}, pp. 767--775.

\bibitem{Odonell}
Khot, S., Kindler, G., Mossel, E., and O'Donnell, R. (2007)
Optimal inapproximability results for MAX-CUT and other two-variable CSPs? {\em SIAM Journal of Computing} {\bf 37} (1), pp. 319-357 (2007).


\bibitem{Natr2001a} {Natterer, F.} (2001) {\em The Mathematics of Computerized Tomography}, SIAM: Society for Industrial and Applied Mathematics, Classics in Applied
    Mathematics.

\bibitem{Peche}
P\'ech\'e, S. (2006) The largest eigenvalues of small rank perturbations of Hermitian random matrices,
{\em Prob. Theo. Rel. Fields} {\bf 134} (1): 127--174 .


\bibitem{Koby2001}
{Rubinstein, J.} and {Wolansky, G.} (2001) Reconstruction of optical surfaces from
  ray data, \emph{Optical Review}, {\bf 8} (4), pp. 281--283.


\bibitem{AmitPNAS} Singer, A. (2008)
\newblock A Remark on Global Positioning from Local Distances.
\newblock {\em Proceedings of the National Academy of Sciences}, {\bf 105} (28):9507-9511.

\bibitem{common-lines}
Singer, A. and Shkolnisky, Y. ({\em submitted for publication}) Three-Dimensional Structure Determination from Common Lines in Cryo-EM by Eigenvectors and Semidefinite Programming.

\bibitem{class-averaging}
Singer, A., Shkolnisky, Y., Hadani, R. ({\em in preparation}) Viewing Angle Classification of Cryo-Electron Microscopy Images using Eigenvectors.

\bibitem{Soshnikov}
Soshnikov, A. (1999)
Universality at the edge of the spectrum in Wigner random matrices.
{\em Comm. Math. Phys.} {\bf 207}, pp. 697--733.

\bibitem{SDPT31}
Toh, K.C., Todd, M.J., and Tutuncu, R.H. (1999)
SDPT3 --- a Matlab software package for semidefinite programming,
{\em Optimization Methods and Software,} {\bf 11}, pp. 545--581.

\bibitem{TracyWidom}
Tracy, C. A., and H. Widom (1994)
Level-spacing distributions and the Airy kernel.
{\em Communications in Mathematical Physics} {\bf 159} (1), pp. 151--174.

\bibitem{SDPT32}
Tutuncu, R.H., Toh, K.C., and Todd, M.J. (2003)
Solving semidefinite-quadratic-linear programs using SDPT3,
{\em Mathematical Programming Ser. B,} {\bf 95}, pp. 189--217.

\bibitem{Boyd1996}
{Vandenberghe, L.}, and {Boyd, S.} (1996)
\newblock Semidefinite programming.
\newblock {\em SIAM Review}, {\bf 38}(1):49--95.

\bibitem{small_world}
Watts, D. J. and Strogatz, S. H. (1998) Collective dynamics of small-world networks. {\em Nature} {\bf 393},
pp. 440--442.

\bibitem{Wigner1}
Wigner, E. P. (1955) Characteristic vectors of bordered matrices with infinite dimensions,
{\em Annals of Mathematics} {\bf 62} pp. 548--564.

\bibitem{Wigner2}
Wigner, E. P. (1958) On the distribution of tile roots of certain symmetric matrices,
{\em Annals of Mathematics} {\bf 67} pp. 325--328.


\bibitem{DC1} Wuthrich, K. (2003). NMR studies of structure and function of biological macromolecules (Nobel Lecture). {\em J Biomol NMR} {\bf 27}, 13-39.






\end{thebibliography}
\end{document}